\documentclass{amsart}
\usepackage{hyperref}
\usepackage[english]{babel} 
\usepackage{graphicx}

\usepackage{amssymb}
\usepackage{amsfonts}
\usepackage{amscd}

\newtheorem{theorem}{Theorem}[section]

\newtheorem{proposition}[theorem]{Proposition}

\newtheorem{definition}[theorem]{Definition}

\theoremstyle{definition}
\theoremstyle{remark}
\numberwithin{equation}{section}

\newcommand{\BB}{{\mathbb B}}
\newcommand{\FF}{{\mathbb F}}
\newcommand{\PP}{{\mathbb P}}
\newcommand{\GG}{{\mathbb G}}
\newcommand{\VV}{{\mathbb V}}
\newcommand{\vol}{\mathbf{vol}}
\newcommand{\dist}{\mathrm{dist}}
\newcommand{\AAA}{{\mathbb A}}

\newcommand{\reals}{\mathbb{R}}

\newcommand{\sgn}{\mathrm{sgn}}
\newcommand{\mlt}{\mathbf{mlt}}

\newcommand{\rationals}{\mathbb{Q}}

\newcommand{\labelpag}[1]{\refstepcounter{enumi}\label{#1}}

\begin{document}
\title[Affine Exterior Algebra]{Reworking on Affine Exterior Algebra \\ of Grassmann: Peano and his School}
\author{Gabriele H. Greco}
\address{Dipartimento di Matematica\\
Universit\`{a} di Trento, 38050 Povo (TN), Italy}
\email{greco@science.unitn.it}
\author{Enrico M. Pagani}
\address{Dipartimento di Matematica\\
Universit\`{a} di Trento, 38050 Povo (TN), Italy}
\email{pagani@science.unitn.it}
\dedicatory{On the occasion of the $150^{th}$ anniversary of the birth of Giuseppe Peano}
\date{%
August 21, 2009
}

\begin{abstract} 
In this paper
a construction
of 
affine exterior algebra of \textsc{Grassmann}, with a special attention to the 
revisitation of this subject operated by \textsc{Peano} and his School,
is 
examined from a historical viewpoint. Even if the
exterior algebra over a vector space is a well known concept, the construction
of an exterior algebra over an affine space, in which points and vectors  
coexist, has been neglected. This paper wants to fill this lack.

Some attention is given to the introduction of defining 
by abstraction (today called definition by quotienting or by equivalence relation),
a procedure due to and used by \textsc{Peano} to define geometric forms, basic
elements of 
an 
affine exterior algebra. This \textsc{Peano}'s innovative 
way of defining, is a relevant contribution to mathematics. 

It is observed that in the construction of an affine exterior algebra on the Euclidean 
three-dimensional space, \textsc{Grassmann} and \textsc{Peano} make use of metric concepts:
an accurate analysis shows that, in some cases, the metric aspects can be eliminated,
putting into evidence the sufficiency of the underlying affine structure of the Euclidean space.

In the final part of the paper some geometrical and mechanical 
applications and interpretations of the affine exterior algebra given by \textsc{Grassmann} 
and \textsc{Peano} are presented.
\end{abstract}

\maketitle

\section{Introduction}

\emph{Calcolo geometrico} \cite{calcolo88} of 1888 marks strong interest of \textsc{Giuseppe Peano} (1858-1932) in \emph{Extension Theory} (\emph{Ausdehnungslehre}) of \textsc{Hermann Grassmann} (1809-1887): despite difficulties related to a philosophical aura~\footnote{Nowadays a philosophical prejudice is nothing but a pretext to keep away from \emph{Ausdehnungslehre}. This prejudice, shared also by some prominent mathematicians, have been originated by the fact that, following the habit of the time, \textsc{Grassmann} received a philosophical-theological education.  \textsc{Enriques} in \cite[(1911) p.\,71]{enriques} says: ``... le syst\`eme philosophique [of \emph{J. F. Herbart}] 
a exerc\'e [...] une grande influence sur le d\'eveloppement des id\'ees de 
\emph{H. Grassmann} et de \emph{B. Riemann}''.
Moreover,  \textsc{Weyl} says in \cite[(1918) p.\,319]{weyl}:
``In forming the conception of a manifold of more than three dimensions, Grassmann as well as Riemann was influenced by the philosophic ideas of Herbart.''.

\textsc{Lewis} \cite[(1977)]{lewis} provides further connections between theology and \emph{Ausdehnungslehre}.
It is worth saying, however, that a comparative analysis of the works of \textsc{Grassmann}, \textsc{Schleiermacher}, \textsc{Herbart} and \textsc{Fichte} shows that this philosophical influence is not perceptible to us.    
} 
of the book of \textsc{Grassmann}, \textsc{Peano} was attracted by the applications of 
the new calculus proposed by \textsc{Grassmann}. 
\footnote{
\textsc{Peano} in \cite[(1896) p.\,953]{saggio} says:
``La teoria di Grassmann \`e oggigiorno, dai varii autori che l'hanno riesposta ed applicata,
giudicata coi pi\`u grandi elogi. [...] se tanto tard\`o quest'opera a farsi conoscere,
e se tanta difficolt\`a presenta tuttora nel diffondersi, la ragione ci dev'essere, e,
secondo me, essa sta nella forma dell'esposizione, forma metafisica e nebulosa, 
lontana dal linguaggio solito dei matematici, e che fin da principio, invece di attirare i lettori,
li stanca ed allontana. Ed anch'io, nello studio di quest'opera, rilevai la potenza del nuovo metodo
solo nell'esame delle applicazioni, specialmente quelle pubblicate nel 1845 
nell' `Archiv' di Grunert ({\em Grassmann's Werke}, I, p.\,297 [English trans. in \cite[p. 283]{grass_branch}])''.  
}

Geometric calculus of \textsc{Grassmann}
consists of an algebraic calculus on geometric entities 
(points, oriented segments, oriented triangles and so on), 
differently from what happens in analytical geometry, that 
is based on an algebraic manipulation of coordinates.  
This distinguished aspect of the geometrical calculus of \textsc{Grassmann} is emphasized by \textsc{Peano} in the introduction of his \emph{Calcolo geometrico}:

\begin{quotation}
Il calcolo geometrico presenta analogie con la geometria analitica; ne differisce in ci\`o, che,
mentre nella geometria analitica i calcoli si fanno sui numeri che determinano gli enti geometrici, in questa nuova scienza i calcoli si fanno sugli enti stessi.
\end{quotation}

Among the members of the School of Peano \cite[p.\,187]{peano_vita}, who share with him the interest for \textsc{Grassmann}, we have to mention \textsc{Burali-Forti}, \textsc{Castellano}, \textsc{Boggio}, \textsc{Bottasso}, 
\textsc{Pensa}, \textsc{Marcolongo} and \textsc{Burgatti}, the last two mainly oriented to mathematical-physics.
All them took an active role in development and diffusion of the ``vector calculus \emph{\`a la Grassmann}''
\footnote{
This active role is also proved by the publication of the series of books 
{\em Analisi vettoriale generale e applicazioni}
\cite[(1930)]{burali_marcolongo_1930}, \cite[(1930)]{burgatti_boggio_burali_1930},
\cite[(1930)]{burgatti_1930}, \cite[(1930)]{boggio_1930}, \cite[(1930)]{marcolongo_1930}
that received a warm approval of \textsc{Peano} \cite[(1915)]{peano_simboli}: 
``I suddetti professori [Burali-Forti e Marcolongo], in unione al Prof.\,Boggio, Bottasso ed altri, hanno cominciato la pubblicazione di una serie di volumi, ove sono trattate le principali applicazioni dei vettori. Sicch\`e i pi\`u bei libri su questa teoria, che una volta erano stampati in Inghilterra, ora sono pubblicati in Italia.''. 
}, hereafter referred to as ``geometric calculus''.

\textsc{Burali-Forti}, in collaboration with \textsc{Marcolongo}, published many papers and books on vector calculus (see for instance
\cite[(1897)]{burali_1897}, \cite[(1909)]{burali_1909}, \cite[(1932)]{burali_1932}). 
The first book of rational mechanics, in terms of concepts introduced by \textsc{Grassmann} and revisited by \textsc{Peano} was written by
\textsc{Castellano} \cite[(1894)]{castellano}.
\textsc{Boggio} wrote some books and papers on vector calculus and differential geometry.
\textsc{Bottasso} must be mentioned for his book \emph{Astatique} \cite[(1915)]{bottasso}.
\textsc{Pensa} extended in \cite[(1919)]{pensa,pensaMI} the results of \emph{Calcolo geometrico} of \textsc{Peano} 
from $3$-dimensional
to arbitrary $n$-dimensional spaces. 

\emph{Calcolo geometrico} of \textsc{Peano}
revisits in an original way the ideas of \textsc{Grassmann};
in particular he fruitfully introduces the modern notion of \emph{vector space}. 
In construction of an exterior algebra
on Euclidean three-dimensional space \textsc{Grassmann} and \textsc{Peano} make use of metric concepts:
an accurate analysis shows that, in a large part of the geometric calculus, the metric attribute of the space can be eliminated, putting into evidence the sufficiency of the underlying affine structure alone.

A concise presentation of the geometric calculus is given by \textsc{Peano} in 
an essay \emph{Elementi di calcolo geometrico}  \cite[(1891)]{peano1891_elementi}
and in a book for university students \emph{Lezioni di analisi infinitesimale} \cite[(1893) pp.\,16-41]{peano1893}; 
a wider comprehensive presentation is given
in an essay \emph{Saggio di calcolo geometrico} \cite[(1895)]{saggio}; 
an axiomatic exposition is given in a paper
\emph{Analisi della teoria dei vettori} \cite[(1897)]{peano_vettori} and in the celebrated
\emph{Formulario mathematico} (5th edition
\cite[(1908) pp.\,188-201, 269-270]{peano_1908}, 4th edition
\cite[(1902-03) pp.\,277-285]{peano_F4}). 
In his seminal book \emph{Applicazioni geometriche} \cite[(1887)]{peano87}, in a paper \emph{Teoremi sui massimi e minimi} \cite[(1888)]{peano_massimi} and in a book for students \emph{Lezioni di analisi infinitesimale} \cite[(1893)]{peano1893}, 
\textsc{Peano} gives 
applications of \textsc{Grassmann}'s exterior algebra.

The aim of this paper, focussed on historical facts 
\footnote{
This paper is not a definitive word on the revisitation of the School of Peano of the 
geometric calculus of \textsc{Grassmann}: only a part of that
will be reworked in this paper. In particular the so called regressive
product of \textsc{Grassmann}, the differential and integral aspects and, in general, all the ``metric'' 
concepts, might deserve a specific paper. Another lack of this paper is the absence of
comparison between  different ways of reception
of \textsc{Grassmann}'s geometric calculus by \textsc{Peano} and  other
authors (for instance, \textsc{Hankel}, \textsc{Schlegel}, \textsc{Whitehead}, \textsc{Cartan}, \textsc{Fehr}, \textsc{Forder}, \textsc{Hyde},~\dots).
Moreover the contributions of the members of the School of Peano to the geometric calculus
should require more attention and further investigations.
We hope to fill these gaps in forthcoming papers.
}, 
is to give a brief presentation of the main ideas
of \emph{affine} aspects of \textsc{Grassmann}'s exterior algebra, a subject that has not 
yet received enough attention from the mathematical community. In our exposition, reorganizing \textsc{Peano}'s contribution,
we offer a reader a logical path in order to make 
the comprehension of geometric calculus easier.

From a methodological
point of view, we are focussed on primary sources, and not on the secondary ones, that is, on mathematical facts, and not on opinions or
interpretations of other scholars of history of mathematics.

In section \ref{affine_vector} we examine concepts of affine and 
vector space and 
of affine volume,
in view of the introduction of an affine exterior algebra.

In Section \ref{def_math_sec}  we present \emph{definitions
by abstraction}, today called definition by quotienting or by equivalence relation.
This  procedure is due to and used by \textsc{Peano} to define geometric forms.
In our opinion, the \textsc{Peano}'s innovative way of definition through an 
equivalence relation, is a relevant contribution to Mathematics. 

In sections \ref{definitions}, \ref{reduction_sec}, \ref{geometric_sec} and
\ref{mechanical_sec}, devoted to affine exterior algebra (definition, properties, 
geometric and mechanical interpretation respectively), we present the 
reconstruction  of the theory of \textsc{Grassmann} due to \textsc{Peano}.

\section{Affine and vector spaces, affine volume and linear extensions of affine spaces} \label{affine_vector}

In \textsc{Grassmann}'s work, points and vectors coexist distinctly in a common
structure, together with other objects, like exterior products of points (for instance, bipoints,
tripoints, k-points) and vectors (for instance, bivectors, trivectors, k-vectors). 
This subtle distinction was very demanding in comparison with today habits of mathematicians.

The conceptual distinction between points and vectors, and consequently between
affine and vector spaces, is a necessary and fundamental
precondition in order to grasp the geometrical ideas of \textsc{Grassmann}. In the second
edition of \emph{Ausdehnungslehre} (1862)
the structure of the exterior algebra, as usual at the present day, is built on vectors; on the contrary, in the first edition of 
\emph{Ausdehnungslehre} (1844) it is built on points. 

\textsc{Peano} maintains firmly the distinction between points and vectors 
and so on. 
In several papers he applies the geometric calculus of \textsc{Grassmann}, for
instance in {\em Applicazioni geometriche} 
(see \cite[(1887) p.\,164]{peano87}) and in \cite[(1890)]{peano_area})
\footnote{
In \emph{Applicazioni geometriche}, published by \textsc{Peano} one year before \emph{Calcolo 
geometrico}, differential and integral calculus for functions with 
values in an \emph{exterior algebra} is developed too.}
he defines the area of a surface, and in \cite[(1898)]{peano_vettori} 
he gives an axiomatic re-foundation (today standard) of affine spaces and 
Euclidean geometry, based on the three primitive notions of \emph{point},
\emph{vector} (i.e., difference of points) and \emph{scalar product}.

In addition to the axiomatization of affine space, \textsc{Peano} in \emph{Calcolo geometrico} 
\cite[(1888), Cap. IX, pag. 141--142]{calcolo88} provides a modern definition of vector space. 

In \textsc{Grassmann} there is no explicit definition of abstract vector space.
Nevertheless the  following properties (\ref{prima})-(\ref{ultima}),  based on \emph{sum}, \emph{difference} and \emph{multiplication by scalar} (see \cite[(1862) pp.\,4--5 and p.\,129]{grass_E}) are recognized  
\emph{fundamental} by \textsc{Grassmann} in order to obtain all the other algebraic properties
depending on the three operations: sum, difference and multiplication. 
It is easy verifying that these properties describe exactly a  vector spaces, as stated in the following 
proposition.

\begin{proposition}
Let $V$ be a nonempty set, endowed with two binary operations ``$+$'', ``$-$'' such that, for every $a,b,c \in V$, the following 
two conditions hold:
\begin{eqnarray}
a+b= b+a,\quad a+(b+c)=(a+b)+c \labelpag{prima}
\\
(a+b)-b =a , \quad 
(a-b)+b=a.              \labelpag{seconda}
\end{eqnarray}
Then the element ${\bf0}:= a-a$ is well defined (it does not depend on the choice of $a\in V$) and  $(V,+)$ is a commutative group whose null element is $\bf 0$. 
\footnote{
In the first edition of \emph{Ausdehnungslehre} \cite[(1844), pp.\,33-40]{grass_branch} 
\textsc{Grassmann} analyzes in details the two operations of addition `$+$' and subtraction `$-$'.
It is worth noticing that properties (\ref{prima})-(\ref{seconda}) characterize
an abelian group.
} 

If, in addition,  $V$  is endowed  with  a left and right multiplication by  real numbers (denoted merely by juxtaposition) such that, for every $a,b\in V$ and $\beta, \gamma$ real numbers, the following three conditions hold:
\begin{eqnarray}
a \beta = \beta a, \quad (a \beta) \gamma = a (\beta \gamma), \quad a 1 = a,\\
(a+b) \gamma = a \gamma + b \gamma, \quad a (\beta + \gamma) = a \beta + a \gamma \\
a \beta = \bf{0}\quad\text{if and only if $a=\bf{0}$ or $\beta = 0$}\,,\labelpag{ultima}
\end{eqnarray}
then  $V$ is a (real) vector space with respect to $+$ and multiplication by scalars.
\end{proposition}

\textsc{Peano}'s definition of vector space is based on the following four primitive ingredients: 
1) \emph{null element}, 2) \emph{equivalence relation}, 
3) \emph{addition} and 4) \emph{multiplication by real numbers}. 

\begin{definition} \label{peano_def_vec}{\rm (\textsc{Peano}'s definition of a vector space \cite[(1888) \S\,72]{calcolo88}).%
\footnote{\textsc{Peano} used the name ``\emph{linear system}''  in place of 
our ``\emph{vector space}''.}} A (nonempty) set $V$, endowed with a binary relation $\cong$, a binary operation $+$
and a left multiplication by real numbers (denoted merely by juxtaposition), is said to be a vector space, if there is an element 
${\bf0}\in V$ (called \emph{null element}) and the following properties hold:
\begin{eqnarray}
a\cong b\iff b\cong a  \label{simmetria} \\
(a\cong b\text{ and } b\cong c)\Longrightarrow a\cong c \label{transitivita} \\
a\cong b\Longrightarrow a+c\cong b+c\\
a+b\cong b+a, a+(b+c)\cong (a+b)+c\\
a\cong b\Longrightarrow m a\cong mb\\
m(a+b)\cong ma+mb, (m+n)a\cong ma+na, m(na)\cong(mn)a,\\ 
1a\cong a, 0a\cong {\bf0}   \label{unita}
\end{eqnarray}
for every $a,b,c \in V$ and $m,n$ real numbers.
\footnote{
It is worth noticing that reflexivity property of the binary relation $\cong$ 
follows by (\ref{simmetria}), (\ref{transitivita}), and (\ref{unita}).
}
\end{definition}

When the equivalence is in fact an identity, \textsc{Peano}'s
definition \ref{peano_def_vec} agrees completely with the modern definition of vector space.
The presence of an equivalence relation, instead of an identity, makes the definition
more powerful, allowing him, for instance, to define a structure of vector spaces 
on a set whose elements are freely generated by fixed rules acting on some given {\em \`a-priori} elements.
 Basic examples are given by the definition of a ``\emph{linear extension}'' of an affine space (see \textsc{Peano} \cite[(1898), pp.\,525-526]{peano_vettori})  and  by the definition of ``\emph{geometric forms}'' in \emph{Calcolo geometrico} \cite[(1888) n.\,5]{calcolo88}.

The totality of geometric forms of \textsc{Peano} is a vector space: its basic four ingredients mentioned
above being defined through the notion of volume of the Euclidean $3$-space. An accurate analysis puts into evidence that the metric aspects of the space are not strictly necessary for a large amount of geometric calculus. 

This ``non-metric'' approach may be pursued 
by introducing a notion of \emph{affine volume}.  
\emph{Affine volume} over an arbitrary $n$-dimensional affine space, is a function 
defined on the ordered $(n+1)$-tuples of vertices of a (possibly degenerate) $n$-simplex,
and it is 
characterized by the following four properties: non-triviality (i.e.\,\,volume of some $n$-simplex is 
non-zero), translational invariance, change of sign under an odd permutation of its vertices, multi-affine dependence on the vertices.

In a modern language this non-metric approach can be formulated in terms of
a \emph{Peano space}, namely 
of a pair $(\AAA_n, \vol_n)$ where $\AAA_n$ is an $n$-dimensional affine space and $\vol_n$ is an
affine volume over $\AAA_n$. As usual, the affine space $\AAA_n$ will be regarded as a triple 
$(\PP_n, \VV_n, +)$, where $\PP_n$ is the underlying set of  \emph{points} of $\AAA_n$, 
$\VV_n$ is the associated  $n$-dimensional \emph{vector space} and $+ : \PP_n \times \VV_n \to \PP_n$
is the usual \emph{translation operation}. 
\footnote{
For given points $p, q \in \PP_n$ the difference $p-q$ denotes the usual difference in an affine space,
namely the unique vector $v \in \VV_n$ such that $p= q+v$.
}

\begin{definition}
The affine volume
${\vol_n}:(\PP_n)^{n+1}\to\reals$ is characterized the following properties:
\begin{enumerate}
\item $\vol_n$ is not vanishing,   \labelpag{vol_1}
\item $\vol_n(p_1+v,\dots,p_{n+1}+v)=\vol_n(p_1,\dots,p_{n+1})$,  \labelpag{vol_2}
\item
$\vol_n(p_{\sigma(1)},\dots,p_{\sigma({n+1})})=\sgn(\sigma) \vol_n(p_1,\dots,p_{n+1})$,  \labelpag{vol_3}
\item$ \vol_n(\alpha p_1+\beta p_1', p_2,\dots,p_{n+1})=
\alpha\vol_n( p_1, p_2,\dots,p_{n+1})+\beta\vol_n(  p_1',     
p_2,\dots , p_{n+1})$    \labelpag{vol_4}
\end{enumerate}
for arbitrary
 $p'_1,p_1,\dots,p_{n+1}\in \PP_n$, $v\in V_n$, $\beta:=1-\alpha$,
$\alpha\in\reals$    and permutations $\sigma$ of $1, ..., n+1$.
\end{definition}

It is worth observing that in the definition of geometric forms, \textsc{Peano}
involves the ``ratios'' between volumes of simplexes: in this way \textsc{Peano}'s definition
is independent of the choice of a specific volume. 
This independence is due to the fact that the
ratio between the volumes of two non-degenerate $n$-simplices is invariant
under affine transformations.

Geometric forms allowed  \textsc{Peano} to give, in different ways, 
\emph{linear extensions} of affine spaces. 
\footnote{
\textsc{Berger} in \cite[(1987), p.\,68]{berger} introduced the notion of \emph{universal vector space}.
\textsc{Berger} says: ``the construction of this universal space may at first appear to come out of the blue''.
On the contrary a linear extension of an affine space (which ``coincides'' with the universal space)
are easily introduced by \textsc{Peano}.
}
Given an affine space $\AAA$ and a vector space $W$,
we say that $W$ is a \emph{linear extension} of $\AAA$ with respect to an injective
application $j: \AAA \to W$ if $j$ is an affine application and its image is a hyperplane $H$ of $W$ with $0 \not\in H$.

The approach of \textsc{Burali-Forti} \cite[(1915)]{burali-forti_1915} to geometric forms
suggests a way of constructing  a linear extension (denoted by $\BB (\AAA_n)$ and called the
\emph{Burali-Forti space} of $\AAA_n$) using the notion of affine volume introduced above.
In this sense, we regard points of $\PP_n$ as
functions on sets of $n$-tuples of points of $\PP_n$ through the application $\vol_n$.
More precisely, for all $p \in \PP_n$, we denote by $j(p) : {(\PP_n)}^n \to \reals$
the function defined by 
$$j(p)(p_1,p_2, ..., p_n) := \vol_n(p, p_1, p_2, ..., p_n)  \, .$$
Then the vector space $\BB (\AAA_n)$ is defined to be the totality of functions $j(p)$. 
By properties (\ref{vol_1} - \ref{vol_4}) the Burali-Forti space $\BB (\AAA_n)$ 
 is a linear extension of the affine space $\AAA_n$ with respect to $j$.

It is not surprising that linear relations in Burali-Forti  space $\BB (\AAA_n)$ can be recovered as 
affine relations in  $\AAA_n$ (and conversely \footnote{
For instance, through a relation of type (\ref{aff_p}) it is possible to characterize
a ``barycenter'': given a (finite) system of weighted points $\{(p_i, \alpha_i)\}$ such that
$\sum_i \alpha_i =1$, its barycenter is the point $G$ uniquely determined by the condition
$\sum_i \alpha_i (p_i - p) = G - p$, for all point $p \in \PP_n$.
Therefore, by (\ref{aff_vol}), the barycenter $G$ may be characterized through the condition 
\begin{equation} \labelpag{Carnot}
\sum_{i=1}^m \alpha_i\,\vol_n(p_i, q_2,\dots,q_{n+1})=\vol_n(G, q_2,\dots,q_{n+1})
\text{ for any } q_2,\dots,q_{n+1}\in \PP_n .
\end{equation}
Condition (\ref{Carnot}) motivates the characterization of the barycenter given by \textsc{Carnot} 
\cite[(1801) p.\,154]{carnot} as ``\emph{centre des moyennes distances}'', i.e., the 
center of the  signed distance of the point $p_i$ (weighted with $\alpha_i$) by 
 arbitrary planes. 
}), in virtue  of the  following proposition.

\begin{proposition} Let $ \{\alpha_i\}_{i=1}^m\subset\reals$ and    \label{equi_bur}
$\{p_i\}_{i=1}^m\subset\PP_n$. The following properties are equivalent:
\begin{enumerate}
\item $\sum_{i=1}^m \alpha_i(p_i-p)=0$ for any  $p\in \PP_n$ \,   \labelpag{aff_p}
\item $\sum_{i=1}^m \alpha_i\,\vol_n(p_i, q_2,\dots,q_{n+1})=0$\quad for any
\quad $q_2,\dots,q_{n+1}\in \PP_n$ .     \footnote{
The proof can be based on the use of the  \emph{vector volume} $\mlt_n : (\VV_n)^n \to \reals$ associated to the affine
volume $\vol_n$ and well defined by $\mlt_n(v_1, ..., v_n) := \vol_n(p, p+v_1, ..., p+v_n)$
for arbitrary $p \in \PP_n$ and $v_1, ..., v_n \in \VV_n$. Observe that the
application $\mlt_n$ is multilinear and alternating. 
}
   \labelpag{aff_vol}
\end{enumerate}
\end{proposition}

Other examples of linear extensions are obtained by introducing 
\emph{M\"obius spaces}.
A pair formed by a vector space W and a non-vanishing linear form $\mu : W \to \reals$
is called a \emph{M\"obius space} with \emph{mass} $\mu$, in honor of \textsc{M\"obius}. 
To a M\"obius space is associated an affine space $\AAA := (\PP,\VV, +)$,
where $\PP := \{x\in W : \mu(x) =1\}$, $\VV := \{ x \in W : \mu(x)=0 \}$ and 
the translation operation $+$ is the restriction to $\PP \times \VV$ of the sum operation 
$+$ on $W$. 
It is clear that a M\"obius space $W$ is a linear extension of the affine space $\AAA$
with respect to the canonical injection from $\PP$ to $W$.

\section{Definitions in Mathematics} \label{def_math_sec}

A relevant subject of the research activity of \textsc{Peano} and his School concerned {\em definitions} in
Mathematics, a subject that received and till now receives more attention by philosophers than 
by mathematicians. 
\footnote{For instance, recently (in the year 2005) the philosophical Journal {\em Mind} has celebrated with a special issue (1222 pages) the Centenary of \textsc{Bertrand Russell}'s
landmark essay {\em On Denoting}  \cite[(1905)]{russell_denoting}.
\textsc{Russell}, attracted by the ingenious formalization of the use
of the article ``the'' in mathematical definitions performed by \textsc{Peano} through the
introduction of the ``\emph{inverse iota descriptor}'', writes the essay {\em On Denoting}
to make precise the use of the article ``the'' from logic and linguistic points of view.
}

The problem of defining has been treated by \textsc{Aristotle} in {\em Analytica Posteriora} (II, 10)
and in {\em Topica} (I, 5).
Among the philosophers who proposed a theory of definition we mention
\textsc{Hobbes} and \textsc{Leibniz} (see \textsc{Couturat} 
\cite[(1901) cap.\,VI, n.\,7]{couturat}).

In 1903 \textsc{Vailati}, a member of the School of Peano,
discovered (see \cite{saccheri})
an interesting booklet: {\em Logica demonstrativa} of \textsc{Saccheri} 
\footnote{
According to  \textsc{Vailati} (\cite[(1903), p.\,326]{aristotele}), the historical and philosophical relevance of the
{\em Logica demonstrativa} of \textsc{Saccheri} is comparable with the \emph{Logique de Port Royal} (\cite[(1662)]{port_royal})
of \textsc{Arnauld} and \textsc{Nicole}.  
}.
Saccheri, in studying the {\em Organon} (or logical works) of \textsc{Aristotle}, recognized
in his \emph{Logica demonstrativa} two modalities of definitions
\cite[(1697) pars II cap. IV]{saccheri_L}: nominal ({\em quid nominis}) and real ({\em quid rei}), giving examples of
mistakes originating by confusing these two types. 
A definition {\em quid rei} denotes an existing object; on the contrary a definition
{\em quid nominis} gives a meaning to a word.
Concerning a definition {\em quid nominis}, its correctness simply relies on the fact that the
word has not been previously defined. In order for the correctness of a definition
{\em quid rei}, \textsc{Saccheri} considers two cases: if the definition is \emph{incomplexa} (i.e., concerning
only one property), it is necessary to postulate or to give a proof of the existence of
{\em definitum}; if the definition is \emph{complexa} (i.e., concerning more properties),
it is necessary to show the compatibility of the properties themselves.

In \cite[(1899)]{vacca} \textsc{Vacca}, another member of the School of Peano, praises the theory of definitions of \textsc{Gergonne} \cite[(1818-19)]{gergonne}, and gives rich comments
 about {\em implicit definitions}. 

\textsc{Peano} is aware of the historical development of definition theory \cite[(1900)]{peano_def1900}. 
In 1894 \textsc{Burali-Forti} in his booklet {\em Logica matematica} \cite{burali-logica} gives a classification
of various types of definitions: 1) nominal definitions without any hypothesis 
(i.e., the {\em definiens} is made only of constant symbols whose meaning
is fixed by previous definitions), 2) nominal definitions
with hypothesis (in which the {\em definiens} contains variables whose meaning
is determined by the hypothesis), 3) definition by induction, 4) implicit definitions and
5) definition by abstraction. 

The term \emph{definition by abstraction} has been introduced in mathematics by \textsc{Peano}
\cite[(1894) \S 38-39]{peano_quotient}. 
Such way of defining, starting from a given set $A$ and a reflexive, symmetric and transitive
relation $\approx$ on $A$, provides a pair ($X, \varphi$) made of a set $X$ and a surjective function 
$\varphi : A \to X$ such that
\begin{equation}   \label{plurality}
\varphi(a) = \varphi(b) \iff a \approx b \, .
\end{equation}

In the review paper \cite[(1915)]{peano_review} \textsc{Peano} gives a rich and significant list,
starting with \textsc{Euclides}, of
mathematical entities constructed through the definition by abstraction (for example, rational numbers, real numbers \footnote{For instance,  let $A$ be the totality of upper bounded subsets of $\rationals$ (the rational numbers). For every pair $a,b\in A$ define $a\approx b$ by ``$a$ and $b$ have the same upper bounds''. The relation $\approx$ is an equivalence relation. By definition by abstraction from $(A, \approx)$, it follows the existence of a pair $(X,\varphi)$.  $X$ ``is'' the set of real numbers and $\varphi$ ``is'' $\sup$ (\emph{supremum}); the algebraic structure  on $X$ is well defined by the \emph{zero} $\textbf{0}:=\varphi(\{0\})$, the \emph{unit} $\textbf{1}:=\varphi(\{0\})$, the \emph{addition} $\varphi(a)+\varphi(b):=\varphi(a+b)$, the \emph{order}  $\big(\varphi(a)\le\varphi(b)\big):=\big(\varphi(a\cup b)=\varphi(b)\big)$, the \emph{product} of positive numbers  $\varphi(a)\cdot\varphi(b):=\varphi(a^{+}\cdot b^{+})$, whenever $a^{+}:=\{q\in a: q>0\}\ne\emptyset$ and  $b^{+}:=\{q\in b: q>0\}\ne\emptyset$.
}, cardinality of sets,  directions, vectors
 and so on).
In the same paper \textsc{Peano} exposes ``alternative ways'' of formulating a definition
by abstraction: in particular he considers the case in which $\varphi$ is defined by 
\begin{equation}  \label{padoa_russell}
\varphi (a):= \text{ ``the set of element } b \in A \text{ such that } b\approx a \text{''}.
\end{equation} 
The definition (\ref{padoa_russell}), ascribed by \textsc{Peano} to \textsc{Russell}, corresponds to the modern concept of {\em quotient space}. \textsc{Russell} introduced  the definition (\ref{padoa_russell})  in \cite[(1901), p.\,320]{russell3} 
(see also \cite[(1903) \S 108-111]{russell_principles});
he observed its relevance in providing nominal definitions whenever definitions by
abstraction are possible, and  in avoiding the plurality of pairs $(X,\varphi)$ involved in 
definition (\ref{plurality}).

Definitions by abstraction have been analyzed from different viewpoints (correctness,
history, avoidance of logical drawbacks, practical use)
by members of the School of Peano:
\textsc{Burali-Forti} \cite[(1899)]{burali-fortiE}, \textsc{Padoa} \cite[(1908)]{padoa}, \textsc{Vailati} \cite[(1908)]{vailati_G},
\textsc{Burali-Forti} \cite[(1912)]{burali-forti_1912}, \textsc{Maccaferri} \cite[(1913)]{maccaferri} and
\textsc{Burali-Forti} \cite[(1915)]{burali-forti_1915},
\cite[(1924)]{burali-forti_1924}, \cite[(1925)]{burali-forti_1925}. 

Definitions by abstraction play a basic role in this paper. 
We will see in the next section that the spaces of geometric forms were
obtained by \textsc{Peano} in \emph{Calcolo geometrico} by quotienting 
a freely generated vector space with respect 
to a suitable subspace. 
This way of defining geometric forms 
 marks the explicit origin of the modality of definition
by abstraction in mathematics. 
According to \textsc{Vailati} \cite[(1908) p.\,105]{vailati_G}, \textsc{Grassmann} with
\emph{Ausdehnungslehre} paved the way for definitions by abstraction.

\section{Affine exterior algebra: definitions and notation} \label{definitions}

Nowadays an exterior algebra is built on a given vector space. 
\textsc{Grassmann}, in the second edition of \emph{Ausdehnungslehre}, built the exterior algebra by imposing an anticommutative \emph{product} 
(called {\em combinatorial product} by him) between  the linear independent generators, called \emph{unit elements}. 

The elements of his exterior algebra are called by \textsc{Grassmann} 
\emph{extensive magnitudes}.
In his construction
the product is not explicitly defined. His procedure, restated in modern terms, 
consists in quotienting with respect to a suitable equivalence relation (accounting for anticommutativity) the associative algebra freely generated by \emph{unit elements}.  

In the first edition of \emph{Ausdehnungslehre} \cite{grass_branch} 
\textsc{Grassmann} built an \emph{affine exterior algebra} from \emph{points}.
A notion of \emph{mass} of elementary magnitudes allows him to 
distinguish points by vectors. Besides, \textsc{Grassmann} introduced  
an operator (called {\em divergence} by \textsc{Grassmann}, denoted by $\omega$ by \textsc{Peano},
and corresponding to the modern {\em boundary operator} $\partial$)
that extends the notion of mass to the whole exterior algebra, enabling him  to separate 
products of points and  products of vectors.

Let $\GG(\AAA_{n})$ denote the affine exterior algebra related to an $n$-dimensional affine  
space $\AAA_{n}$. Following \textsc{Grassmann}, in $\GG(\AAA_{n}$) there are elements of zero degree (i.e., real numbers), 
elements of degree 1 (i.e., points of $\AAA_{n}$ and their linear combinations), elements of degree 2
(i.e., products of two points of $\AAA_{n}$ and their linear combinations) and so on.

For $k=1,2, \dots $,  the symbol ${\mathbb F}_k(\AAA_{n})$ denotes the linear 
combinations of products of $k$ points of $\AAA_{n}$, called {\em geometric forms of degree} $k$ by 
\textsc{Peano}.
${\mathbb F}_0$ denotes the set of real numbers.

The operator $\omega$ is defined on all geometric forms. The linearity and  the following equalities  characterize $\omega$ uniquely: 
for every points $P_{0}, P_{1}, \dots, P_{k},\dots$ of the affine space $\AAA_{n}$ 
\begin{equation}\label{operator_omega}
\begin{array}{rcl}
\omega(1) &=& 0 \\
\omega(P_{0}) &=& 1  \\
\omega(P_{0}P_{1}) &=& P_{1}-P_{0} \\
\omega(P_{0}P_{1} P_{2}) &=& (P_{1}-P_{0})(P_{2}-P_{0}) \\
\omega(P_{0}P_{1} P_{2}P_{3}) &=& (P_{1}-P_{0})(P_{2}-P_{0})(P_{3}-P_{0}) \\
\dots\\
\omega(P_{0}P_{1}P_{2} \cdots P_{k}) &=& (P_{1}-P_{0})(P_{2}-P_{0})\cdots(P_{k}-P_{0})
\end{array}
\end{equation}

Let $\omega_r$ denote  the \emph{restriction} of $\omega$ to $\FF_r(\AAA_{n})$. 
The operator $\omega_{1}$, the so-called \emph{mass},  allows us to recover 
in  ${\mathbb F}_1(\AAA_{n})$ an affine copy of the affine space 
$\AAA_{n}$ and to deduce that $\FF_{1}(\AAA_{n})$ has dimension $n+1$.
%
In fact,  $\FF_1 (\AAA_n)$ is a linear extension of the affine space $\AAA_n$ with
respect to the canonical injection; in other words, the pair $(\FF_1(\AAA_n), \omega_1)$  is a M\"obius space 
with $\PP_n = \{ x \in  \FF_1 (\AAA_n) : \omega_1(x) =1 \}$.
\footnote{
In the sequel we will see that the linear structure of $\FF_1(\AAA_{n})$
allows also to recover the \emph{barycentric calculus} of \textsc{M\"obius} \cite[(1827)]{mobius}. 
}

Now denote by ${\mathbb V}_k (\AAA_{n})$, $k=1,2, \dots $ the linear 
combinations of products of $k$ vectors of $\GG(\AAA_{n})$, called the {\em space of $k$-vectors}. ${\mathbb V}_0(\AAA_{n})$ denotes the set of real numbers.
 
Clearly, from ($\ref{operator_omega}$) it follows that the operators $\omega_{r}$ maps the geometric forms of degree 
$r$ onto $\VV_{r-1} (\AAA_{n})$.
 
 The elements of $\PP_n$ (resp.\,of  $\VV_{1}(\AAA_{n})$)  are called \emph{points} (resp.\,\emph{vectors}) of the affine exterior algebra $\GG(\AAA_{n})$;  the latter may be  ``identified'' with the vectors of the affine space $\AAA_{n}$. 
 
Following \textsc{Peano}, the product of two, three and four vectors of $\GG(\AAA_{n})$ are called \emph{bivector}, \emph{trivector} and \emph{quadri-vector}, respectively; similarly, the product of two, three and four points of $\GG(\AAA_{n})$ are called \emph{bipoint}, 
\emph{tripoint} and \emph{quadri-point}, respectively.   
Moreover, for any geometric form $x$ of degree $r=2, 3, 4$,   the value $\omega_{r}(x)$ is called by \textsc{Peano} a \emph{vector} (resp. \emph{bivector}, \emph{trivector}) of the geometric form $x$. 

The exterior product is anticommutative, namely 
\[xy = (-1)^{rs} yx \quad \text{for every} \,\, x \in \FF_r(\AAA_n),  y \in \FF_s(\AAA_n) \,.
\]
In particular, the product of two geometric forms of degree $1$ is null if the factors are 
linearly dependent.
Therefore ${\VV}_k(\AAA_{n})$ contains only the zero geometric form, for $k > n$. Analogously, ${\FF}_k(\AAA_{n})$ contains only the zero geometric form, for $k > n+1$; in symbols
\begin{equation}\GG(\AAA_{n})=\FF_{0}(\AAA_{n})\oplus\FF_{1}(\AAA_{n})\oplus\cdots\oplus\FF_{n}(\AAA_{n})\oplus\FF_{n+1}(\AAA_{n}).
\end{equation}

The main innovation operated by \textsc{Peano} in this context can be found in the
construction of an  affine exterior algebra  based on the notion of {\em affine volume}. Concerning $\GG(\AAA_{3})$, \textsc{Peano} says in \emph{Calcolo geometrico} \cite[(1888) p.\,VI]{calcolo88}:

\begin{quotation} 
Le definizioni introdotte per le formazioni di quarta specie [$\FF_{4}(\AAA_{3})$], o volumi, sono gi\`a comuni in geometria analitica; le definizioni per le formazioni delle tre prime specie [$\FF_{1}(\AAA_{3}), \FF_{2}(\AAA_{3}), \FF_{4}(\AAA_{3})$] sono ridotte con metodo uniforme a quelle date pei volumi [\dots].

\end{quotation}

Let $\vol_n$ be an affine volume on $\AAA_{n}$. Following \textsc{Peano}, for $1 \le k \le n+1$, the space $\FF_k(\AAA_{n})$ of geometric forms  of degree $k$ is defined by quotienting
the vector space of the \emph{homogeneous polynomials} of the same degree $k$ (freely generated by points)
by the subspace of homogeneous polynomials  
\begin{equation}  \label{polynomial}
\sum_{i} \alpha_{i} P_{(i,1)} \cdots P_{(i,k)}
\end{equation}
such that,  for any choice of the points $P_{k+1}, P_{k+2},\dots , P_{n+1}$ \footnote{For $k=n+1$,  the quantified variables $P_{k+1}, \dots , P_{n+1}$ disappear in formula $(\ref{volume})$.}, the following equality holds true:
\begin{equation}   \label{volume}
\textstyle\sum_{i} \alpha_{i} \vol_n (P_{(i,1)}, \dots , P_{(i,k)}, P_{k+1}, \dots , P_{n+1}) =0.
\end{equation}

To complete the definition of affine exterior algebra given by \textsc{Peano} notice that:
\begin{enumerate}
\item
the exterior product of forms of degree $k$ by forms of degree $s$ (when $0 \le k+s \le n+1$) is defined by \textsc{Peano} starting from the formal product of polynomials,
verifying that the definition is compatible with the operation of quotienting
defined above, and satisfies the
properties of associativity, distributivity and anticommutativity.
\item \textsc{Peano}, following \textsc{Grassmann}, indicates the exterior product by juxtaposition.
\footnote{
In the various editions of \emph{Formulario Mathematico} \textsc{Peano} uses the letter
$\alpha$ to denote the exterior product; for example $P \alpha Q$, $P \alpha Q \alpha R$ in place of  $PQ$
 and $PQR$ respectively. 
}
 Notice that this
notation does not create any ambiguity with the ``polynomial'' notation for 
geometric forms (\ref{polynomial}), because a ``monomial'' geometric form 
$P_{(i,1)} \cdots P_{(i,k)}$ is equal to the exterior product of the points 
$P_{(i,1)}, \dots, P_{(i,k)}$.
\end{enumerate}

The notion of volume, that stands at the basis of the construction of the affine exterior 
algebra, leads quite naturally \textsc{Peano} to describe linear forms  on the space $\FF_k(\AAA_{n})$
of the geometric forms of degree $k$ by geometric forms
of supplementary degree $s:=(n+1)-k$ (see \cite[(1888)  n.\,82]{calcolo88}). 
For instance, given a geometric form 
$\varphi := \sum_{i} \alpha_{i} Q_{(i,1)} \cdots Q_{(i,s)}$ of degree $s$,  a linear 
form $\varphi^{*}$ on $\FF_k(\AAA_{n})$  can be uniquely defined through the relation
\begin{equation}  \label{volume_b}
\varphi^{*} (P_{1}, \dots , P_{k}) := 
\sum_{i} \alpha_{i} \vol_n (P_{1}, \dots , P_{k}, Q_{(i,1)}, \dots ,Q_{(i,s)})
\end{equation}
for any point $P_1, \dots , P_{k}$.
\footnote{It is surprising to perceive in (\ref{volume_b}) a trace of
the idea of Hopf algebra and Hodge dual. 
In this context, it is also worth mentioning the \emph{metric} Hodge $*$ operator that corresponds to the \textsc{Grassmann}'s {\em index}.
}

It is worth observing that \textsc{Burali-Forti} \cite[(1915)]{burali-forti_1915} defines the elements (points, bipoints and so on) that generate the exterior algebra, 
as functions on elements of supplementary degree. More precisely he uses formula
(\ref{volume_b}) to define these functions. 

\textsc{Burali-Forti}'s approach avoids the necessity of introducing the spaces of geometric forms
as quotients (as \textsc{Peano} did in his construction).
Indeed, in \textsc{Burali-Forti} the geometric forms are regarded as functions, thus inheriting the
usual vector space structure. See the Burali-Forti spaces of Section 2.

\textsc{Peano}'s approach to exterior algebra through the notion of volume is emphasized
in section 35 ``\emph{Begr\"undung der Punktrechnung durch G. Peano}''
of \cite[(1923), pp.\,1543-1545]{lotze_1923} in the celebrated 
\emph{Encyklop\"adie der mathematischen Wissenschaften}.
  
In several papers \textsc{Peano} introduces the geometric forms of degree $1$ directly,
without using the notion of volume $\vol_n$. See, for instance, \textsc{Peano} \cite[(1898), p.\,525-526]{peano_vettori}, 
\cite[(1899), p.\,155]{peano_F2}, \cite[(1901), p.\,196]{peano_F3}.
More clearly, \textsc{Peano} describes the space $\FF_1(\AAA_n)$   by quotienting 
 the vector space freely generated by points of $\AAA_n$ with respect to
 the subspace of the free linear combinations $\sum_i \alpha_i P_i$ such that
 \begin{equation}
 \sum_i \alpha_i (P_i - P) = 0 \text{ for all } P \in \PP_n  \,\, .
 \end{equation}

The equivalence between this description of the geometric forms of degree $1$ 
and the previous one given in terms of volume is guaranteed by Proposition \ref{equi_bur}.

\section{Affine exterior algebra $\GG(\AAA_{3})$: properties and reduction formulae}
\label{reduction_sec}

In accordance with the presentation given by \textsc{Grassmann} in 
\emph{Ausdehnungslehre} and by \textsc{Peano} in \emph{Calcolo geometrico},
we restrict our presentation to the affine exterior algebra  $\GG(\AAA_{3})$,  related to a $3$-dimensional affine  
space $\AAA_{3}$:
\begin{equation}\GG(\AAA_{3})=\FF_{0}(\AAA_{3})\oplus\FF_{1}(\AAA_{3})\oplus\FF_{2}(\AAA_{3})\oplus\FF_{4}(\AAA_{3}) \, \,.
\end{equation}
As planned in the Introduction, in this section (and in the subsequent ones) we will
present affine exterior algebra as revisited by \textsc{Peano}, focussing our attention on the
aspects and properties that, in our opinion, are characteristic and reveal the richness, fecundity and 
clarity of the geometric calculus. 
For convenience of the reader some significant proofs of
the simplicity of geometric calculus and some comments are confined in footnotes 
\footnote{
In some cases we have indicated some concepts present in the works of \textsc{Peano}
with a new terminology (for instance, applied vector, applied bivector, applied trivector,
Poinsot pair, boundary-cycle).  
}.

\smallskip{\sc Bases, dimension and coordinates} (see \textsc{Peano} \cite[(1888) \S 60]{calcolo88}, \cite[(1908) p.\,194-195]{peano_1908}). Tetrahedrons suggest a way to construct 
bases for the forms of different degrees
of the affine exterior algebra $\GG(\AAA_{3})$. 

Let fix a (non-degenerate) tetrahedron with vertices $A,B,C,D$.
Its $4$ vertices are a basis of $\FF_1(\AAA_{3})$.  The 6 bipoints $AB, AC, AD, BC, BD, CD$ which correspond to its 6 edges, 
 are a basis of $\FF_2(\AAA_{3})$. To its faces there correspond 4 tripoints
 $ABC, ACD, ABD, BCD$ which are a basis of $\FF_3(\AAA_{3})$. Finally, 
 the quadri-point $ABCD$  is a basis of the  $1$-dimensional $\FF_4(\AAA_{3})$. 
\footnote{Similarly,
in an $n$-dimensional affine space $\AAA_{n}$ one has \[\dim\big(\FF_k(\AAA_{n})\big)={n+1\choose k},\] for $0 \le k \le n+1$. If  $T$   is an arbitrary $n$-simplex of $\AAA_{n}$, a basis of $\FF_k(\AAA_{n})$  corresponds to  ${n+1\choose k}$  $k$-dimensional facets of $T$.}

The coordinates of a geometric form of degree $1$ with respect to 
the basis of $\FF_{1}(\AAA_{3})$
formed by the four vertices of a tetrahedron are called {\em baricentric coordinates}. The coordinates  with respect to an arbitrary basis  are referred  to as \emph{projective coordinates}.
Choosing as a basis of $\FF_{1}(\AAA_{3})$ the point $A$ and the
$3$ vectors $B-A$, $C-A$, $D-A$, one has a (non-orthogonal) \emph{Cartesian
coordinate system} with origin $A$ and axes directed along the three vectors.

 Let  $x_{1}$, $x_{2}$, $x_{3}$ and $x_{4}\in \FF_{1}(\AAA_{3})$ be arbitrary geometric forms of degree 1. They are a basis of $\FF_{1}(\AAA_{3})$ if and only if $x_{1}x_{2}x_{3}x_{4}\ne 0$. 
 
 Now assume that  $x_{1}$, $x_{2}$, $x_{3}$, $x_{4}  \in \FF_{1}(\AAA_{3})$ is a basis of  $\FF_{1}(\AAA_{3})$.

 \begin{enumerate}  
\item
Then the product $x_{1}x_{2}x_{3}x_{4}$ is a basis of the $1$-dimensional space $\FF_{4}(\AAA_{3})$. 
 For all $z \in \FF_{4}(\AAA_{3})$ let us denote by
 \begin{equation*}
 \frac{ z}{x_{1}x_{2}x_{3}x_{4}}.
 \end{equation*}
 the real number $\lambda$ such that $z = \lambda \, x_{1}x_{2}x_{3}x_{4}$.
 
 \item If $\alpha_{1}$, $\alpha_{2}$, $\alpha_{3}$ and $\alpha_{4}$ are the coordinates of an element $ x \in \FF_{1}(\AAA_{3})$
 with respect to the basis $x_{1}$, $x_{2}$, $x_{3}$ and $x_{4}$, then the following formulae
\footnote{These formulae are obtained  multiplying the two members of the equality $ x=\alpha_{1}x_{1}+\alpha_{2}x_{2}+\alpha_{3}x_{3}+\alpha_{4}x_{4}$ by the tripoints $x_{2}x_{3}x_{4}$, $x_{1}x_{3}x_{4}$, $x_{1}x_{2}x_{4}$, and $x_{1}x_{2}x_{3}$, respectively. A similar calculus provides coordinate formulae  (\ref{coord_2}) and (\ref{coord_3}).} hold true
 \begin{equation*}
 \alpha_{1}=\frac{ xx_{2}x_{3}x_{1}}{x_{1}x_{2}x_{3}x_{4}}, \, \,\alpha_{2}=\frac{x_{1} xx_{3}x_{4}}{x_{1}x_{2}x_{3}x_{4}}, \, \,\alpha_{3}=\frac{x_{1}x_{2} xx_{4}}{x_{1}x_{2}x_{3}x_{4}}, \, \, \alpha_{4}=\frac{x_{1}x_{2}x_{3} x}{x_{1}x_{2}x_{3}x_{4}}.
 \end{equation*}

 \item Moreover, $x_{1}x_{2}$, $x_{1}x_{3}$, $x_{1}x_{4}$, $x_{2}x_{3}$, $x_{2}x_{4}$, $x_{3}x_{4}$ is  a basis of $\FF_{2}(\AAA_{3})$; the corresponding coordinates $\alpha_{12}, \alpha_{13}, \alpha_{14}, \alpha_{23}, \alpha_{24}, \alpha_{34}$ of an element 
 $s  \in \FF_{2}(\AAA_{3})$ satisfy the following formulae \labelpag{coord_2}
 
 \begin{eqnarray*}
 \alpha_{12}=\frac{ sx_{3}x_{4}}{x_{1}x_{2}x_{3}x_{4}}, \,\, \alpha_{13}=\frac{ sx_{4}x_{2}}{x_{1}x_{2}x_{3}x_{4}}, \,\ \alpha_{14}=\frac{sx_{2}x_{3}}{x_{1}x_{2}x_{3}x_{4}},\\ \alpha_{23}=\frac{ sx_{1}x_{4}}{x_{1}x_{2}x_{3}x_{4}},
\,\, \alpha_{24}=\frac{ sx_{3}x_{1}}{x_{1}x_{2}x_{3}x_{4}},\,\,\alpha_{34}=\frac{ sx_{1}x_{2}}{x_{1}x_{2}x_{3}x_{4}}.
 \end{eqnarray*}

\item The set of geometric forms $x_{2}x_{3}x_{4}$, $-x_{1}x_{3}x_{4}$, $x_{1}x_{2}x_{4}$ and $-x_{1}x_{2}x_{3}$ is  a basis of $\FF_{3}(\AAA_{3})$; the corresponding coordinates $\beta_{1}, \beta_{2}, \beta_{3}, \beta_{4}$ of an element 
$y \in \FF_{3}(\AAA_{3})$ satisfy the following formulae \labelpag{coord_3} 
 
 \begin{equation*}
 \beta_{1}=\frac{x_{1} y}{x_{1}x_{2}x_{3}x_{4}}, \,\, \beta_{2}=\frac{x_{2} y}{x_{1}x_{2}x_{3}x_{4}},\,\, \beta_{3}=\frac{x_{3} y}{x_{1}x_{2}x_{3}x_{4}}, \,\,\beta_{4}=\frac{x_{4} y}{x_{1}x_{2}x_{3}x_{4}}. 
 \end{equation*}

 \end{enumerate}

\smallskip{\sc The operator $\omega$ and the formula of reduction}. 
In addition to linearity, properties of  the operator $\omega$ 
on products of geometric forms of any degree are contained in the following exhaustive list \cite[(1908) p.\,197]{peano_1908}: 
\footnote{\label{foot__Leibniz}
Nowadays it is relevant the notion of {\em differential graded algebra}. Such an algebra
is endowed with a linear operator $\partial$ satisfying the following properties:
\begin{equation*} 
\begin{array}{llcl}
(a) &\partial(\partial(x))  &=& 0   \\
(b) &\partial(xy) &=& \partial(x)y +(-1)^r x \partial(y)  \\
\end{array}
\end{equation*}
for any $x, y$ belonging to the algebra, with $x$ of degree $r$.

It is worth observing that last two properties are satisfied by the operator $\omega$.
The first one (a) is merely a consequence of (\ref{im_ker}). 
We found no explicit evidence of the second property (b) in the works of \textsc{Grassmann} and \textsc{Peano}, but it is evident that properties (\ref{Peano_Leibniz})
are nothing but instances of (b).
}
\begin{equation} \label{Peano_Leibniz}
\begin{array}{lcl}
\omega(x_{1}x_{2}) &=& \omega(x_{1}) x_{2} - x_{1} \omega (x_{2}),  \\
\omega(x_{1}x_{2}x_{3}) &=& \omega(x_{1}) x_{2}x_{3} + \omega(x_{2})x_{3}x_{1} + \omega(x_{3}) x_{1}x_{2}
       , \\
\omega(xs) &=& \omega(sx) = \omega(x)s + \omega(s)x 
        \\     
\omega(xy) &=& \omega(x) y - x \omega(y), 
       \\
\omega(s_{1}s_{2}) &=& \omega(s_{1})s_{2} + s_{1} \omega(s_{2}) ,
      \\
\omega(x_{1}x_{2}x_{3}x_{4}) &=& \omega(x_{1})x_{2}x_{3}x_{4} - \omega(x_{2})x_{1}x_{3}x_{4} + \omega(x_{3})x_{1}x_{2}x_{4} - \omega(x_{4})x_{1}x_{2}x_{3}
\end{array}
\end{equation}
for every $x, x_{1}, x_{2}, x_{3}, x_{4}\in \FF_{1}(\AAA_{3})$, $s, s_{1}, s_{2}\in\FF_{2}(\AAA_{3})$, $y\in\FF_{3}(\AAA_{3})$. 

Concerning the problem of reducing geometric forms, the following fundamental 
\emph{reduction formula} holds:
\begin{equation}\label{red_formula}
x=P\omega(x)+\omega(Px)
\end{equation}
for every point $P$ and every geometric form $x$, whose proof follows immediately
from (\ref{Peano_Leibniz}) or, alternatively, by the linearity of $\omega$.
\footnote{
Instances of this formula, 
that encompasses all reduction formulae present in \textsc{Peano}, will be detailed in the following.
We note that reduction formula is a straightforward consequence of the 
signed Leibniz rule (b) of footnote \ref{foot__Leibniz}.}

 A consequence of the reduction formula is
\begin{equation}
\VV_{r}(\AAA_{3}) = {\rm im} \,\omega_{r+1}=\ker \, \omega_{r}     \label{im_ker}.
\end{equation}

\begin{figure}[htbp]
\begin{center}
\[
\begin{array}{l|ccccccccccccc}
\text{dimension} &0&&1&&4&&6&&4&&1&&0\cr
\hline
\text{geometric forms} &0&\stackrel{\omega_{5}}{\to}&\FF_{4}&\stackrel{\omega_{4}}{\to}&\FF_{3}&\stackrel{\omega_{3}}{\to}&\FF_{2}&\stackrel{\omega_{2}}{\to}&\FF_{1}&\stackrel{\omega_{1}}{\to}&\FF_{0}&\stackrel{\omega_{0}}{\to}&{0}\cr
\hline
\text{inclusions}&&&\uparrow&&\uparrow&&\uparrow&&\uparrow&&\uparrow&\cr
\VV_{r}={\rm im} \,\omega_{r+1}=\ker \, \omega_{r}&0&&0&&\VV_{3}&&\VV_{2}&&\VV_{1}&&\VV_{0}&&0\cr
\text{dimension}&&&0&&1&&3&&3&&1&&0\cr\end{array}
\]
\caption{Geometric forms in a 3-dimensional affine space}
\label{tabella}
\end{center}
\end{figure}

The equality (\ref{im_ker}) summarizes the three basic properties of operator $\omega$:

\begin{enumerate}
\item  ``$\VV_{r}(\AAA_{3}) = {\rm im}\,\omega_{r+1}$'', that is:  $\omega_{1}$ sends
(weighted) points into scalars (their weight), $\omega_{2}$ sends bipoints into vectors, $\omega_{3}$ send tripoints
into bivectors, and  $\omega_{4}$ send quadri-points into trivectors.     

\item ``$\VV_{r}(\AAA_{3}) \subset\ker \,\omega_{r}$'', that is: the mass of a vector, the vector of a bivector,  the bivector of a trivector are all null.
\item ``$\VV_{r}(\AAA_{3}) \supset\ker \,\omega_{r}$'', that is: a geometric form $x$ of degree $r=1$ (resp. $2$ and  $3$) is a vector (risp. bivector and trivector), whenever  its mass $\omega_{1}(x)$ (resp. vector $\omega_{2}(x)$ and bivector $\omega_{3}(x)$) is null.

\end{enumerate}

\smallskip

{\sc Geometric forms of degree $1$} (see \textsc{Peano} \cite[(1888) chap.\,II]{calcolo88}, \cite[(1908) p.\,196]{peano_1908}, \cite[(1893) nn.\,265-266]{peano1893}).
An instance of the reduction formula $(\ref{red_formula})$ is the following reduction formula for a form $x$ of degree $1$:
\begin{equation}\label{red_1}
x  = P \omega_1 (x) +\omega_2(Px) \quad\text{for every point } P.
\end{equation} 
Consequently, for every  point $P$, a form $x \in \FF_1(\AAA_{3})$ can be reduced to a sum of a weighted point 
$P \omega_1 (x)$ and a vector $\omega_{2}(Px)$. If $x=\sum_{i}\alpha_{i}P_{i}$,  the equality $(\ref{red_1})$ becomes:
\begin{equation}
\sum_{i}\alpha_{i}P_{i}= (\sum_{i}\alpha_{i})P+\sum_{i}\alpha_{i}(P_{i}-P).
\end{equation}
If $x\ne0$ and  the  vector  $x - P \omega_1 (x)$, called {\em deviation} by \textsc{Grassmann}, is null, then  the point $P$ is said to be a \emph{barycenter} of $x$.  

 A geometric form $x$ of degree $1$ is a vector, if its mass $\omega_{1}(x)$ is null; in such a case, for every point $P$, the element $Q:=x+P$ is a point  and
\begin{equation}
x=Q-P;
\end{equation}
otherwise, if $\omega_1(x)\ne 0$,  $x$  is a weighted point, that is
\begin{equation}
x =\omega_{1}(x)\big(\frac{x}{\omega_{1}(x)}\big).
\end{equation}
and the element $\frac{x}{\omega(x)}$ is a point: the barycenter of $x$. \smallskip

{\sc Geometric forms of degree $2$}. (see \textsc{Peano} \cite[(1888) chap.\,III]{calcolo88}, \cite[(1908) p.\,197]{peano_1908}, \cite[(1893) nn.\,267-270]{peano1893}).
An instance of the reduction formula $(\ref{red_formula})$ is the following reduction formula for a form $x$ of degree $2$: 
\begin{equation}  \label{red_2}
x =  P \omega_2 (x) + \omega_3(Px).
\end{equation} 
Consequently, a form $x \in \FF_2(\AAA_{3})$ can be reduced to a sum of a bipoint 
$P \omega_2 (x)$ and a bivector $\omega_3(Px)$.
If $x=\sum_{i}Q_{i}P_{i}$, then the equality $(\ref{red_2})$ becomes:
\begin{equation}
\sum_{i}Q_{i}P_{i}=P \sum_{i}(P_{i}-Q_{i})+\sum_{i}(Q_{i}-P)(P_{i}-P).
\end{equation}

There are different useful ways of representing bipoints and bivectors. Name \emph{applied vector} every product $Av$ of a point $A$ and a vector $v$.  The difference  
\[Av-Bv\] is said to be a \emph{Poinsot pair}, whenever $v$ is a vector and  $A, B$ are two points.  A geometric form $x$ of degree $2$ is said to be a \emph{boundary-cycle of a triangle} \footnote{A boundary-cycle of a triangle can be viewed as oriented boundary of a triangle.}
, if there exists three points $A,B,C$
such that
\begin{equation}\label{ciclo}
x=AB+BC+CA
\end{equation}

Moreover
\begin{enumerate}
\item
The bipoints are the {applied vectors}.\footnote{For every bipoint $AB$ it is clear that $B-A$ is a vector and $AB=A(B-A)$. Conversely,
let $A$ be a point and $v$ a vector, then $Av = AB$, where $B:=A+v$.}
\item
The bivectors are the {Poinsot pairs}.\footnote{If $v$, $w$ are vectors and $A$ is
an arbitrary point, then $v w = B w - A w$ where $B := A+v$. Conversely,  a Poinsot pair  $Av-Bv$  may be rewritten as $(A-B)v$.}
\item
The bivectors are the {boundary-cycles of triangles}.\footnote{If $v$, $w$ are vectors and $A$  is
an arbitrary point, then $v w = (B -A)(C-A)=AB+BC+CA$ where $B := A+v$ and $C:=A+w)$. Conversely,  a boundary-cycle $AB+BC+CA$  may be rewritten as $(B -A)(C-A)$.}
\item Every couple of bivectors can be reduced to a couple of bivectors with a common factor. \footnote{That is: if  $v_{1}, v_{2}, w_{1}, w_{2}$ are vectors, then there exist vectors $y,v,w$ such that $v_{1}v_{2}=yv$ and $w_{1}w_{2}=yw$. }
\item
A linear combination of bivectors is a bivector ({\em Poinsot rule}).

\item A geometric form $x\in \FF_2(\AAA_{3})$  is a bivector if and only if $\omega_{2}(x)=0$.\labelpag{bi_carat}
 \item The vector space $\VV_2(\AAA_{3})$  has dimension $3$.
Given a basis
$v_1, v_2, v_3$ of $\VV_1(\AAA_{3})$, the exterior products $v_1 v_2$, $v_2 v_3$, $v_3 v_1$, 
form a basis for $\VV_2(\AAA_{3})$.

\end{enumerate}

    From the reduction formula $(\ref{red_2})$ and the description $(\ref{ciclo})$ of a bivector, it follows that, for every geometric form $x$ of degree $2$, there exist four points $P,B,C,D$ such that
\begin{equation}x=PB+PC+CD+DP.
\end{equation}
Hence every geometric form of degree $2$ is one of the following three types:
\begin{equation}
\begin{cases}
\text{(i) a bipoint:}& PB\\
\text{(ii) a bivector:}&PC+CD+DP\\
\text{(iii) neither  bipoint nor  bivector:}&{PB+PC+CD+DP \text{ with }P,B,C,D\atop
\text{vertices of a non-degenerate tetrahedron}}
\end{cases}
\end{equation}

In particular, 
every geometric form of degree $2$ is either a bipoint or a bivector, if and only if 
\begin{equation}
xx=0.
\end{equation}

\smallskip

{\sc Geometric forms of degree $3$} (see \textsc{Peano} \cite[(1888) chap.\,IV]{calcolo88}, \cite[(1908) p.\,197]{peano_1908}, \cite[(1893) nn.\,271]{peano1893}).
Another instance of the reduction formula $(\ref{red_formula})$ is the following reduction formula for a form $x$ of degree $3$: 
\begin{equation}  \label{red_3}
x =  P \omega_3 (x) + \omega_4(Px).
\end{equation} 
Consequently, a form $x \in \FF_3(\AAA_{3})$ can be reduced to a sum of a tripoint 
$P \omega_2 (x)$ and a trivector $\omega_4(Px)$.
If $x=\sum_{i}R_{i}Q_{i}P_{i}$, then the equality $(\ref{red_3})$ becomes:
\begin{equation}
\sum_{i}\alpha_{i}Q_{i}P_{i}R_{i}= \sum_{i}P(P_{i}-Q_{i})(R_{i}-Q_{i})+\sum_{i}(Q_{i}-P)(P_{i}-P)(R_{i}-P).
\end{equation}

Name \emph{applied bivector} every product $Ay$ of a point $A$ and a bivector $y$.  The difference  $Ay-By$ is said to be a \emph{Poinsot triangle pair}, whenever $y$ is a bivector and  $A, B$ are two points.  A geometric form $x$ of degree $3$ is said to be a \emph{boundary-cycle of a tetrahedron} \footnote{A boundary-cycle of a tetrahedron can be viewed as oriented boundary of a tetrahedron.}, if there exist 
four points $A,B,C, D$
such that
\begin{equation}\label{ciclo2}
x=BCD - ACD + ABD - ABC
\end{equation}

\begin{enumerate}
\item

The tripoints are the {applied bivectors}.\footnote{For every tripoint $ABC$ it is clear that $(B-A)(C-A)$ is a bivector and $ABC=A(B-A)(C-A)$. Conversely,
let $A$ be a point and $v, w$  vectors, then $Avw = ABC$, where $B:=A+v$ and $C:=A+w$.}
\item
The trivectors are the {Poinsot triangle pairs}.\footnote{To prove this property use the following equality: $(B-A) wz = B wz - A wz$.}
\item
The trivectors are the {boundary-cycles of tetrahedron}.\footnote{To prove this property use the following equality: $(B -A)(C-A)(D-A)=BCD - ACD + ABD - ABC$.}

\item A geometric form $x\in \FF_3(\AAA_{3})$  is a trivector if and only if $\omega_{3}(x)=0$. \labelpag{tri_carat}

\item The vector space $\VV_3(\AAA_{3})$  has dimension $1$.
Given a basis
$v_1, v_2, v_3$ of $\VV_1(\AAA_{3})$, the exterior product $v_1 v_2v_3$ 
is a basis for $\VV_3(\AAA_{3})$. 

\item A geometric form $x \in \FF_3(\AAA_{3})$ is either a tripoint or a trivector.

\end{enumerate}

\smallskip

{\sc Geometric forms of degree $4$}.
Any form $x$ of degree $4$ is a quadri-point, that may be regarded as a product
of an arbitrary point $P$ by a trivector:
\begin{equation}   \label{red_4}
x= P \omega_4 (x) 
\end{equation}

Quadri-points are related to volume forms
and allow to check linear independence of forms $x_1, x_2, x_3, x_4$ of degree $1$:
\begin{equation}
x_1 x_2 x_3 x_4 \ne 0  \Longleftrightarrow \, \, x_1, x_2, x_3, x_4 \, 
\text{are linearly independent}.
\end{equation}
 
Due to anticommutativity of the exterior product, the vector space $\VV_{4}(\AAA_{3})$ is $0$-dimensional, that is 
\begin{enumerate}
\item the exterior product of four vectors is null.
\end{enumerate}
This dimensionality restriction can be restated in terms of points: 
for every points $P, A, B,C,D$, the following equality holds:
\begin{enumerate}
\item $ABCD-PBCD+PACD-PABD+PABC=0$,\label{dimensionalita}\footnote{The left side of equality (\ref{dimensionalita}) can be viewed as oriented boundary of a 4-simplex.}
\end{enumerate}
since $(A-P)(B-P)(C-P)(D-P)=ABCD-PBCD+PACD-PABD+PABC$.

\smallskip{\sc Relevant examples of reductions} (see \textsc{Peano} \cite[(1888) nn.\,27, 31]{calcolo88}, \cite[(1887) n.19, 24]{peano87}, \cite[(1893) nn.\,269, 273]{peano1893}).
In various moments of his activity, for the purpose of evaluating areas and volumes, \textsc{Peano} considered reduction formulae for (plane or not) polygons and for (open or closed) polyhedral surfaces. Two significant examples of reductions,
due to the characterizations of bivectors and trivectors (see ($\ref{bi_carat}$) and 
($\ref{tri_carat}$), respectively) are:

\begin{enumerate}
\item   \labelpag{polygons}
Let $ s:= \sum_{i} A_iB_i $ denote a geometric form of degree $2$ by $A_i,B_i$ points of a plane $\pi$ and with $\omega(s)=0$.
There exist $X,Y,Z$ in $\pi$ such that    
\[\sum_i PA_iB_i = XYZ
\]    
for any choice of the point $P$ in $\pi$.
\footnote{For the proof, the condition $\omega(s)=0$ implies that $s$ is a bivector;
therefore there exist three points $X,Y,Z$ in $\pi$ such that $s= (Y-X)(Z-X)$.
Multiplying by $P$ we get the equality $Ps= P(Y-X)(Z-X)$;  moreover, being $P(Y-X)(Z-X)=XYZ$
(because $X,Y,Z,P$ belong to the same plane $\pi$) and $Ps=\sum_i PA_iB_i$,
we have the required equality.
}

\item  \labelpag{polyhedral}
Let $x:= \sum_{i=1}^n A_iB_iC_i$ denote a geometric form of order $3$ such that 
$\omega (x) =0$. 
There exist four points $X,Y,Z,T$ such that 
\[
\sum_i PA_iB_i C_i= XYZT
\]      
for any choice of the point $P$.
\footnote{For the proof, the condition $\omega(x)=0$ implies that $x$ is a trivector;
therefore there exist four points $X,Y,Z,T$ such that $x= (Y-X)(Z-X)(T-X)$.
Multiplying by $P$ we get the equality $Px= P(Y-X)(Z-X)(T-X)$;  moreover, being 
$P(Y-X)(Z-X)(T-X)=XYZT$ and $Px=\sum_i PA_iB_i C_i$,
we have the required equality.
}
\end{enumerate}

\section{Affine exterior algebra $\GG(\AAA_{3})$: geometrical interpretation}
\label{geometric_sec}

Accordingly to the following citation of \textsc{Peano}, the nature of the relation between algebraic and geometric objects is made clear by examining the equalities between geometric forms:
\begin{quotation}
Il calcolo geometrico, in generale, consiste in un sistema di operazioni a eseguirsi
su enti geometrici, analoghe a quelle che l'algebra fa sopra i numeri. Esso permette di esprimere con formule i risultati di costruzioni geometriche, di rappresentare con equazioni proposizioni di geometria, e di sostituire una trasformazione di equazioni ad un ragionamento.
\end{quotation}

Following \textsc{Grassmann}, \textsc{Peano} associates a same geometric concept to different elements of the affine exterior algebra: this is why different attributes (length, area, volume, direction, orientation,
position and so on) of a geometric object may be taken into account in the various cases.
For instance, \textsc{Peano} associates to points, bipoints, tripoints, quadri-points, the geometric objects: points, segments, tri-angles, tetrahedra respectively; similarly to vectors, bivectors, trivectors, he associates segments, triangles, tetrahedra. 

In what follows, by segment, triangle, tetrahedron we will denote the geometric sets
of three-dimensional affine space $\AAA_{3}$
convexly generated by their vertices (two, three, four respectively). 
Moreover, with oriented segment, triangle, tetrahedron we will denote the pair formed by
the geometric sets with an order of their vertices.
We will denote an oriented geometric object of this type through the ordered list of its vertices.
It is evident that the notion of affine volume given above in section \ref{affine_vector} concerns oriented tetrahedra.

For a better understanding of the point of view of \textsc{Peano} 
(that is, in a pre-vectorial epoch) 
it is necessary to regard the usual affine space $\AAA_{3}$ as a set of points endowed with the action of \emph{translations}. The notions of 
parallelism between two straight lines, between a straight line and a plane and
between two planes, direction of straight lines and of planes) and affine transformations
\footnote{In a pre-vectorial epoch were regarded as \emph{collinearities}.}
are well defined in this context.

Concerning orientation, \textsc{Peano} introduces the following definitions:

\begin{enumerate}

\item
two (non-degenerate) oriented tetrahedra have the \emph{same orientation} if the affine volumes of the corresponding tetrahedra have the same sign.

\item
two (non-degenerate) oriented triangles of vertices $A,B,C$ and $Q,R,S$ belonging to the same plane $\pi$ have the \emph{same orientation} if, for any point $P \notin \pi$, the tetrahedra 
of vertices $A,B,C,P$ and $Q,R,S,P$ have the same orientation.

\item
two (non-degenerate) oriented segments of vertices $A,B$ and $Q,R$ belonging to the same straight line 
$r$ have the \emph{same orientation} if, for any points $P,S \notin r$, the tetrahedra 
of vertices $A,B,P,S$ and $Q,R,P,S$ have the same orientation.
\end{enumerate}

Concerning extension (i.e., length, area, volume), \textsc{Peano} retains the metric definitions that are standard in
Euclidean geometry (choice of a unit measure, orthogonality and so on). A deeper 
analysis of \textsc{Peano}'s geometric calculus reveals that his definitions relies only on  purely
affine properties of the  underlying space. 
Consequently, we will assume, as a starting point, the following purely affine definitions, that are consistent with the use of the metric notion of extension given by \textsc{Peano}:

\begin{enumerate}

\item
two oriented tetraheda have the \emph{same extension} if their affine volumes have the same absolute value.

\item
two oriented triangles of vertices $A,B,C$ and $Q,R,S$ belonging to the same plane $\pi$ have the \emph{same extension} if, for any point $P \notin \pi$, the oriented tetrahedra 
$A,B,C,P$ and $Q,R,S,P$ have the same extension.

\item
two oriented segments of vertices $A,B$ and $Q,R$ belonging to the same straight line $r$ have the \emph{same extension} if, for any point $P,S \notin r$, the oriented tetrahedra 
$A,B,P,S$ and $Q,R,P,S$ have the same extension.
\end{enumerate}

Following \textsc{Peano}, equalities between elements of the affine exterior algebra $\GG(\AAA_{3})$ may be easily interpreted geometrically in $\AAA_{3}$,
as it is evident by the definitions of geometric forms through the notion of affine volume.
More explicitly, we have the following facts:

\begin{enumerate}

\item
$\sum_{i} \alpha_{i} A_{i} = \sum_{j} \beta_{j} Q_{j}$ if and only if
for arbitrary points $R,S,T$ 
the sums of the affine volumes
of the two systems of oriented tetrahedra 
$A_{i},  R, S,T$ and $Q_{j}, R,S,T$ coincide, \labelpag{momento-statico}
 that is
\begin{equation*}
\sum_{i} \alpha_{i}  \, A_{i} R S T=\sum_{j} \beta_{j} \, Q_{j}R S T  \, \, .
\end{equation*}

\item
$\sum_{i} \alpha_{i} A_{i}B_{i} =\sum_{j} \beta_{j} Q_{j}R_{j}$
if and only if for arbitrary points $S,T$ the weighted sums of the affine volumes
of the two systems of oriented tetrahedra 
$A_{i},B_{i},S,$ $T$ and $Q_{j},R_{j},S,T$ coincide.

\item
$\sum_{i} \alpha_{i} A_{i}B_{i}C_{i}=\sum_{j} \beta_{j} Q_{j}R_{j}S_{j} $
if and only if for arbitrary point $T$ the weighted sums of the affine volumes
of the two systems of oriented tetrahedra 
$A_{i},B_{i},C_{i},$  $T$ and $Q_{j},R_{j},S_{j},T$ coincide.

\item
$\sum_{i} \alpha_{i} A_{i}B_{i}C_{i}D_{i}=\sum_{j} \beta_{j} Q_{j}R_{j}S_{j}T_{j}$
if and only if the weighted sums of the affine volumes
of the two systems of oriented tetrahedra 
$A_{i},B_{i},C_{i},D_{i}$ and $Q_{j},R_{j},S_{j},T_{j}$ coincide.

\end{enumerate}

As stated by \textsc{Peano}, the correspondence between elements in the affine exterior algebra $\GG(\AAA_{3})$
and geometric concepts in the usual affine tri-space $\AAA_3$ 
is summarized in the following list of equivalences, of direct proof
using the properties and definitions listed above.

\begin{enumerate}


\item Two points $A$ and $B$ coincide if and only if
$AB = 0$.

\item Three points $A,B,C$ are collinear in $\AAA_3$ if and only if 
$ABC=0$.

\item Four points $A,B,C,D$ are co-planar in $\AAA_3$ if and only if 
$ABCD =0$.

\item
 Two oriented segments of vertices $A,B$ and $C,D$ 
lie on a same straight line and have a same extension and orientation if and only if $AB = CD$.

\item
Two oriented segments of vertices $A,B$ and $C,D$ 
lie on parallel straight line and have a same extension and orientation if and only if $B-A=D-C$.
\footnote{It is worth observing that
\begin{equation} \label{equidifference}
B-A=D-C  \Longleftrightarrow
(B+C)/2 = (D+A)/2 
\end{equation}
This shows that the quadrilateral of vertices $A,B,D,C$
is a parallelogram because the two diagonals meet in their middle point.
The property (\ref{equidifference}) was and is frequently assumed as a starting point of the definition of 
vector as an equivalence class (see \cite[(1898)]{peano_vettori}).
}

\item
Two (non-degenerate) oriented triangles of vertices $A,B,C$ and $P,Q,R$
lie on the same plane and have the same extension and orientation if and only if $ABC = PQR$.

\item
Two (non-degenerate) oriented tetrahedra of vertices $A,B,C,D$ and $P,Q,R,S$
have the same extension and orientation  if and only if $ABCD = PQRS$.

\end{enumerate}

The geometric meaning of the equalities between  vectors, bivector, and trivectors  is deduced from the corresponding equalities  between  bipoint, and tripoints and quadri-points since, for arbitrary vectors $v,w,u,t, r, z$ and a point $P$, the following properties hold:
 \begin{itemize}
\item
$v=w$ if and only if $Pv=Pw$

\item
$vw = ut$ if and only $Pvw=Put$

\item
$vwz = utr$ if and only if $Pvwz=Putr$.

\end{itemize}  

Several identities may be straightforwardly generated by algebraic manipulations. In the
following we list some of them, geometrically significant, present in \emph{Calcolo geometrico} of \textsc{Peano}: 

\begin{enumerate}

\item
A straight line determined by two different points $A$ and $B$  
is parallel to the plane determined by a non-null tripoint $\pi$ if and only if
\[
(B-A) \pi =0.
\]

\item
A straight line determined by two different points $A$ and $B$ is
parallel to the straight line determined by a non-null bipoint $a$ if and only if
\[
(B-A) a =0.
\]

\item
 For every point $P$ belonging to the straight line determined by two different points $A$ and $B$, the following equality holds
\footnote{The equality follows by reduction formula (\ref{red_2}):  
$AB=P \omega(AB) + \omega(PAB) = P(B-A) + \omega (PAB) 
=  P(B-A)$ being $PAB=0$ due to collinearity.}
\[AB=P(B-A) \, .
\]

\item
For every point $P$ belonging to the plane determined by three non-collinear points $A,B,C$, the following equality holds
\footnote{The equality follows by reduction formula (\ref{red_3}):  
$ABC=P \omega(ABC) + \omega(PABC) = P(B-A)(C-A) + \omega (PABC) 
=  P(B-A)(C-A)$ being $PABC=0$ due to co-planarity.}
\[
ABC = P(B-A)(C-A) \, .
\]

\item
$ABCD = P(B-A)(C-A)(D-A)$ for every point $P$ in $\AAA_3$.
\footnote{The equality follows by reduction formula (\ref{red_4}):  
$ABCD=P \omega(ABCD)  = P(B-A)(C-A)(D-A)$.}

\end{enumerate}

Expanding the products in the equalities listed above we have the following description of point belonging to a straight line and to a plane (see \textsc{Peano} \cite{calcolo88}):

\begin{enumerate}
\item Let $A,B$ be two different points and $P$ an arbitrary point. Then 
\[
AB + BP + PA = 0
\]
  if and only if  $P$ belongs to the straight line determined by $A, B$.
\item Let $A,B, C$ be three non-collinear points and $P$ an arbitrary point. Then
\[
ABC + BAP + CBP + ACP  =0 
\]
  if and only if  $P$ belongs to the plane determined by $A,B,C$.
\end{enumerate}

In previous section, using the reduction formula, we presented two statements  (\ref{polygons}) and
(\ref{polyhedral}), which concern polygons  and polyhedral surfaces, respectively.

The geometric interpretation of the first statement corresponds to the following proposition:

\begin{enumerate}
\item \labelpag{pol-bel}
Let $\{A_i \}_{i=1}^{n+1} $ denote a set of points belonging to the plane $\pi$, with
$A_1 = A_{n+1}$. The sum of the areas of the oriented triangles of vertices
$P, A_i, A_{i+1}$, for $1 \le i \le n$, does not depend on the choice of the point
$P$ in the given plane $\pi$.
\end{enumerate}

With respect to the second statement (\ref{polyhedral}) the corresponding geometric interpretation
is:

\begin{enumerate}
\item
Let $\{A_i B_i C_i\}_{i=1}^{n} $ denote a set of tripoints belonging to the affine tridimensional space, with
$\sum_{i=1}^n (B_i - A_i)(C_i - A_i) =0$. 
\footnote{
Whenever 
the set of oriented triangles $A_i B_i C_i$ for $1 \le i \le n$ forms an oriented closed 
polyhedral surface, the condition $\sum_{i=1}^n (B_i - A_i)(C_i - A_i) =0$ is satisfied,
since  $\sum_{i=1}^n (B_i - A_i)(C_i - A_i) =   \sum_{i=1}^n (A_i B_i + B_i C_i + C_i A_i)$. 
}
The sum of the volumes of the oriented tetrahedra of vertices
$P, A_i, B_i, C_i$, for $1 \le i \le n$, does not depend on the choice of the point
$P$ in the given tri-space.
\end{enumerate}

As noticed in \cite{greco-mazzucchi}, these two properties have been proved with slight differences 
by \textsc{Bellavitis} and \textsc{M\"obius}.

\section{Affine exterior algebra $\GG(\AAA_{3})$: mechanical interpretation}
\label{mechanical_sec}

In this section we present \textsc{Peano}'s point of view on elements of the exterior algebra $\GG(\AAA_{3})$  in terms of mechanical concepts.

{\sc Systems of material points and barycenters.} (See \textsc{Peano} \cite[(1893), v. II \S 263]{peano1893}).

Geometric forms of degree $1$ are related to systems $S$ of material points.

Let $\{(P_{i},\alpha_{i})\}$ denote a set of material points $P_i$ of mass $\alpha_i$.
The notion of {\em static moment with respect to a given plane} $\pi$ of a system of points is well known in Physics and is given by $\sum \alpha_i \dist (P_i, \pi)$, sum of the
weighted signed distance of points $P_i$ from $\pi$.
It is evident that this static moment is proportional to the quantity
\begin{equation} 
\sum \alpha_i \vol_n(P_i, A,B,C)
\end{equation}
for arbitrary (not collinear, of course) points $A, B, C$ belonging to $\pi$.

Moreover two material systems of points are said to be (mechanically) \emph{equivalent}
if their \emph{static moments} with respect to any given plane are the same.
A material system having static moment null with respect to any plane,
is called a {\em null system}.

In accordance with the definition of forms of degree $1$ given by \textsc{Peano}
(see \ref{momento-statico}),
we have that 
\begin{quotation}
\emph{two systems of material points are equivalent
if and only if the corresponding geometric forms of degree $1$ coincide.}
\end{quotation}
 In particular to a zero geometric form
of degree $1$ corresponds a \emph{null system}. 

On the basis of the definition of equivalence of material systems of points,
we can observe that a geometric form of degree $1$ with a non-vanishing mass
$\alpha$ coincides with a weighted point $G$ of mass $\alpha$; therefore
$G$ is the physical barycenter. 

{\sc Systems of (applied) forces and equilibrium.}
Following \textsc{Grassmann},  \textsc{Peano} relates  geometric forms $s:= \sum_{i} \alpha_{i} P_i Q_i$ of degree 
$2$ to  systems of applied forces
\footnote{More precisely to the bipoint $ \alpha_{i}P_{i} Q_{i}$ it is associated the force $F_{i}$ 
given by the vector $ \alpha_{i}(Q_{i} - P_{i})$ applied in  $P_{i}$,
and, consequently, to a geometric form $\sum_{i} \alpha_{i} P_i Q_i$, the system of applied forces
$\{F_{i} \}_{i}$ .}.

In Mechanics the \emph{equivalence} of two systems of forces relies on the equality of the
so-called \emph{axial moments} of the two systems with respect to an arbitrary axis.
It is evident that the axial moment is proportional
\footnote{It is worth noticing that the effective \emph{value} of axial moment is a scalar
which depends on the choice
of a metric structure of the space.}
to the quantity
\begin{equation} \label{axial_mom}
\sum \alpha_i \vol_n(P_i, Q_{i},A,B)
\end{equation}
where $A,B$ are two distinct points of the axis. Therefore
\begin{quotation}
\emph{two systems of applied forces are (mechanically) equivalent if and only if the
corresponding geometric forms of degree $2$  coincide.}
\end{quotation}
This fact follows easily by the
definition of geometric form.

As a consequence of this correspondence between 
systems of applied forces and geometric forms, \textsc{Grassmann} and \textsc{Peano} derive some well known statements concerning equivalence of systems of applied forces, reducibility and equilibria:

\begin{enumerate}

\item
a system of applied forces is equivalent (in infinite ways) to  a combination of a force and a \emph{Poinsot pair}
of forces (\emph{Poinsot theorem});  \labelpag{qui} 
 \footnote{ 
Statement (\ref{qui}) is a consequence of the reduction formula (\ref{red_2}) for geometric forms.
Here and in the following ``Poinsot pair of forces'' means ``pair of forces with null resultant''.
The concept of ``pair of forces with null resultant'' is due to \textsc{Poinsot}.

It is worth noticing that geometric result (\ref{pol-bel}) was proved in a \emph{mechanical way} by \textsc{Bellavitis}
(see \cite[(1834) p.\,262]{Bellavitis}) using the notion of \textsc{Poinsot}'s pair of forces.}

\item
the sum of finitely many \emph{Poinsot pair} of forces is a \emph{Poinsot pair};
\footnote{See \emph{Poinsot rule}.}

\item
a planar system of applied forces is equivalent to a force or to a \emph{Poinsot pair} of forces;

\item
a system of parallel applied forces is equivalent to a single force or to a \emph{Poinsot pair} of forces;
\footnote{
In the case of gravitational forces the system of applied forces reduces to a single force applied in the center of mass.}

\item
a system of applied forces is equivalent to a force or to a \emph{Poinsot pair} of forces
if and only if the \emph{scalar invariant} is null;  \footnote{
Denoting with $s$ the geometric form of degree $2$ corresponding to a given system of applied forces,
by reduction formula (\ref{red_2}) there exist a bipoint $a$ and a bivector $\iota$
such that $s= a+ \iota$. Observe that $ss = (a+\iota)(a+\iota) = 2 a \iota$
and $a \iota$ is proportional to the well known \emph{scalar invariant} of the system of applied forces;
therefore $ss=0$ if and only if $a=0$ or $\iota =0$.}

\item
a system of applied forces is reducible to the sum of two forces;

\item
a system of applied forces is reducible in a unique way to a system of forces
acting along the edges of a given (non-degenerate) tetrahedron;
\footnote{This fact follows easily by observing that
the space $\FF_2$ of geometric forms of degree $2$ is $6$-dimensional and
a basis is provided by the $6$ bipoints corresponding to the $6$ edges of a given
tetrahedron.}

\item
characterization of the barycenter $G$ of a material system of points $\{(P_i, \alpha_i)\}_i$
of total mass $\sum_{i} \alpha_i =1$,
through the introduction of  a system of parallel forces (*) or
a system of concurrent forces (**):
\[\begin{array}{lll}
\sum_i \alpha_i P_i = G  &\Longleftrightarrow \quad (*) & \sum_i \alpha_i P_i u =G u \text{ for any vector } u \\\cr
   &\Longleftrightarrow \quad (**) &  \sum_i \alpha_i O P_i  =OG  \text{ for any point } O. 
\end{array} \]
\end{enumerate}

It is surprising that the statements listed above concerning equivalence, reducibility of systems
of applied forces (and in general in problems of static equilibrium), do not require any metric concept 
(i.e., scalar product and vectorial moment). The object of \textsc{Grassmann}-\textsc{Peano} calculus
corresponding to the ``modern'' vectorial moment
is the \emph{Poinsot pair}. 
In \textsc{Grassmann}-\textsc{Peano} calculus applied forces and their moments are
seen \emph{homogeneously}: they are represented by bipoints and Poinsot pairs, namely both of them as forms of degree $2$.
Unfortunately there is no trace of this attitude in the modern presentation of the theory of applied forces:
forces and moments are (at least for their physical dimensions) inhomogeneous objects.

\section{Appendix}

All articles of \textsc{Peano} are collected in \emph{Opera Omnia} \cite{peano_omnia}, a CD-ROM,  edited by \textsc{Roero}. Selected works of \textsc{Peano} were assembled and commented in \emph{Opere scelte} \cite{peano_opere} by \textsc{Cassina}, a disciple of \textsc{Peano}. A few have English translations in \emph{Selected Works} \cite{peano_english}. Regrettably, even fewer \textsc{Peano}'s articles have a public URL and are freely downloadable. 

One finds the following articles of \textsc{Peano}
 \bigskip
 
 \noindent\begin{tabular}{ll}
 \hline
 in \emph{Opere scelte}, vol. 1:& \begin{minipage}{.65\textwidth}
 
  \cite[(1890)]{peano_area}  
 \end{minipage}\cr
\hline

in \emph{Opere scelte}, vol. 2: & \begin{minipage}{.65\textwidth} \cite[(1894)]{peano_quotient}, \cite[(1900)]{peano_def1900}, \cite[(1915)]{peano_review}, 

\cite[(1917)]{peano_eguale}, \cite[(1921)]{peano_def1921},
\end{minipage}\cr
\hline

in \emph{Opere scelte}, vol. 3: & \begin{minipage}{.65\textwidth} \cite[(1888)]{peano_massimi}, \cite[(1891)]{peano1891_elementi}, \cite[(1894)]{peano_recensione_Grassmann},

\cite[(1895)]{peano_recensione_Castellano},
\cite[(1896)]{saggio}, \cite[(1898)]{peano_vettori}, \cite[(1915)]{peano_simboli}
\end{minipage} \cr
\hline

 in \emph{Selected Works}: & \begin{minipage}{.65\textwidth}\cite[(1887) pp.\,152--160, 185--7]{peano87}, \cite[(1887) pp.\,1--32]{calcolo88}, 
 
 \cite[(1890)]{peano_area}, 
 \cite[(1896)]{saggio}, \cite[(1915)]{peano_simboli},
 
 \cite[(1915)]{peano_def1921}. \end{minipage}\cr
\hline
 
 \end{tabular}

\bigskip

For reader's convenience, we provide a chronological list of some
mathematicians mentioned in the paper, together with biographical sources.

The \texttt{html} file  with biographies of mathematicians listed below with an asterisk can be attained at University of St
Andrews's web-page 

\texttt{http://www-history.mcs.st-and.ac.uk/history/\{Name\}.html}
\medskip

\textsc{Arnauld}, Antoine (1612-1694) (*)

\textsc{Nicole}, Pierre (1625-1695)

\textsc{Huygens}, Christiaan (1629-1695) (*) 

\textsc{Leibniz}, Gottfried Wilhelm (1646-1716) (*) 

\textsc{Saccheri}, Girolamo Giovanni (1667-1733) (*) 

\textsc{Carnot}, Lazare Nicolas Margu\'erite (1753-1823) (*) 

\textsc{Gergonne}, Joseph Diaz (1771-1859) (*) 

\textsc{Poinsot}, Louis (1777-1859) (*)

\textsc{M\"obius}, August Ferdinand (1790-1868) (*)

\textsc{Bellavitis}, Giusto (1803-1880) (*)

\textsc{Grassmann}, Hermann (1809-1877) (*) 

\textsc{Hankel}, Hermann (1839-1873) (*)

\textsc{Schlegel}, Victor (1843-1905) see \textsc{May} \cite{may} and \emph{Ens.\,Math.}\,\textbf{8} (1906), p.\,55

\textsc{Hyde}, Edward Wyllys (1843-1930), see \href{http://math.uc.edu/pdf/rtangle03.pdf}{\texttt{http://math.uc.edu/pdf/rtangle03.pdf}}

\textsc{Klein}, Felix (1849-1925) (*)

\textsc{Peano}, Giuseppe (1858-1932), see  \textsc{kennedy} \cite{peano_vita} 

\textsc{Castellano}, Filiberto (1860-1919), see  \textsc{kennedy} \cite{peano_vita}

\textsc{Whitehead}, Alfred North (1861-1947)  (*)

\textsc{Burali-Forti}, Cesare (1861-1931) (*)

\textsc{Marcolongo}, Roberto (1862-1943),  see  \textsc{kennedy} \cite{peano_vita}

\textsc{Vailati}, Giovanni (1863-1909)  (*)

\textsc{Couturat}, Louis (1868-1914) (*)

\textsc{Burgatti}, Pietro (1868-1938),  see \textsc{May} \cite{may}

\textsc{Cartan}, \`Elie (1869-1951) (*)

\textsc{Maccaferri}, Eugenio (1870-1953),  see  \textsc{kennedy} \cite{peano_vita}

\textsc{Fehr}, Henri (1870-1954), see \textsc{May} \cite{may}

\textsc{Enriques}, Federigo (1871-1946) (*)

\textsc{Russell}, Bertrand Arthur William (1872-1970) (*)

\textsc{Pensa}, Angelo (1875-1960), see \textsc{Sallent} \cite{sallent}, p.\,24

\textsc{Vacca}, Giovanni (1875-1953) (*) 

\textsc{Boggio}, Tommaso (1877-1963) (*)

\textsc{Bottasso}, Matteo (1878-1918) (*)

\textsc{Weyl}, Hermann (1885-1955) (*)

\textsc{Forder}, Henry George  (1889-1981), see \emph{Newsletter NZ Math.\,Soc.} \textbf{27} (1983)

\textsc{Cassina}, Ugo (1897-1964), see  \textsc{kennedy} \cite{peano_vita}


\begin{thebibliography}{999}

\bibitem{port_royal} A.~Arnauld, P.~Nicole: 
\newblock {\em Logique de Port Royal},
\newblock Paris, 1662.

\newblock \url{http://gallica.bnf.fr/ark:/12148/bpt6k574432}

\bibitem{Bellavitis} G.~Bellavitis:  \newblock Teoremi generali per determinare le aree dei poligoni e dei volumi e dei poliedri col mezzo delle distanze dei loro vertici,
\newblock {\em Ann. Scienze Regno Lombardo-Veneto}, \textbf{4} (1834), pp.\,256-264.

\newblock \url{http://books.google.it/books/download/Annali_delle_scienze_del_regno_Lombardo_.pdf?id=nuMTAAAAQAAJ&output=pdf&sig=ACfU3U0HjaBity5nT-_zAXkCTF5uMz7abQ}

\bibitem{berger} M.~Berger: 
\newblock \emph{Geometry}, Vol. I,
\newblock Springer-Verlag, Berlin, 1987.

\bibitem{boggio_1930} T.~Boggio: 
\newblock \emph{Analisi vettoriale generale e applicazioni}, Vol. IV {\em Idrodinamica},
\newblock Zanichelli, Bologna, 1930.

\bibitem{bottasso} M.~Bottasso: \newblock {\em Astatique},
\newblock Pavia, 1915.


\newblock
\url{http://quod.lib.umich.edu/cgi/t/text/text-idx?c=umhistmath;cc=umhistmath;view=toc;idno=ABR4766.0001.001}

\bibitem{burali-logica} C.~Burali-Forti: \newblock {\em Logica matematica},
\newblock Hoepli, Milano, 1894, (2nd edition, 1919).

\newblock\url{http://www.archive.org/details/logicamatematic00buragoog}

\bibitem{burali_1897} C.~Burali-Forti: \newblock {\em Introduction \`a la g\'eom\'etrie diff\'erentielle suivant le m\'ethode de H. Grassmann},
\newblock Gauthier-Villars, Paris, 1897.

\newblock \url{http://historical.library.cornell.edu/cgi-bin/cul.math/docviewer?did=Bura001&seq=1}

\bibitem{burali-fortiE} C.~Burali-Forti: \newblock Sur l'\'egalit\'e et sur l'introduction des \'el\'ements  d\'eriv\'es dans la science, 
\newblock {\em Enseignement Math.}, \textbf{1} (1899), pp.\,246--261.

\newblock \url{http://retro.seals.ch/cntmng?type=pdf&rid=ensmat-001:1899:1::88&subp=hires}

\bibitem{burali_1909} C.~Burali-Forti, R.~Marcolongo: \newblock {\em Elementi di calcolo Vettoriale con
numerose applicazioni alla geometria, alla meccanica e alla fisica-matematica},
\newblock Zanichelli, Bologna, 1909 (2nd edition, 1921)

\newblock\url{http://www.archive.org/details/elementidicalcol00burauoft}

\bibitem{burali-forti_1912} C.~Burali-Forti: \newblock Gli enti astratti definiti come enti relativi ad un campo di nozioni, 
\newblock {\em Rend. Acc. Lincei Sc. Fis. Mat. Nat.}, \textbf{21} (1912), pp.\,677--682.

\bibitem{burali-forti_1915} C.~Burali-Forti:  \newblock Nuove applicazioni degli operatori,
\newblock {\em Atti R. Acc. Torino}, \textbf{50} (1915), pp.\,669--684.

\bibitem{burali-forti_1924} C.~Burali-Forti:  \newblock A proposito dell'articolo di E.\,Maccaferri,
\newblock {\em Boll. di Matematica}, \textbf{20} (1924), pp.\,128--129.

\bibitem{burali-forti_1925} C.~Burali-Forti:  \newblock A proposito di una lettera di Mario Pieri,
\newblock {\em Boll. di Matematica}, \textbf{21} (1925), pp.\,136--137.

\bibitem{burali_marcolongo_1930} C.~Burali-Forti, R.~Marcolongo: 
\newblock \emph{Analisi vettoriale generale e applicazioni}, Vol. I {\em Trasfomazioni lineari},
\newblock Zanichelli, Bologna, 1930.


\bibitem{burali_1932} C.~Burali-Forti: \newblock {Elementi di calcolo vettoriale},
\newblock in {\em Enciclopedia delle matematiche elementari e complementi},
Hoepli, Milano, 1932, Vol.\,II--2, pp.\,105--139.

\bibitem{burgatti_boggio_burali_1930} P.~Burgatti, T.~Boggio, C.~Burali-Forti: 
\newblock \emph{Analisi vettoriale generale e applicazioni}, Vol. II {\em Elementi di geometria differenziale} (linee, superficie, spazi curvi, geometria proiettiva differenziale),
\newblock Zanichelli, Bologna, 1930.

\bibitem{burgatti_1930} P.~Burgatti: 
\newblock \emph{Analisi vettoriale generale e applicazioni}, Vol. III {\em Teoria matematica dell'elasticit\`a},
\newblock Zanichelli, Bologna, 1930.

\bibitem{carnot} L.~Carnot: \newblock {\em De la corr\'elation des figures de g\'eom\'etrie},
\newblock Paris, 1801.

\newblock
\url{http://www.archive.org/details/delacorrlationd00carngoog}

\bibitem{castellano} F.~Castellano: \newblock {\em Lezioni di meccanica razionale},
\newblock Torino, 1894 ($2^{nd}$ edition 1911).

\bibitem{couturat} L.~Couturat: \newblock {\em La logique de Leibniz d'apr\`es
des documents in\`edits},
\newblock Paris, 1901 (reprinted by Olms, Hildesheim, 1985).

\newblock \url{http://gallica.bnf.fr/ark:/12148/bpt6k110843d}

\bibitem{enriques} F. Enriques: \newblock {Principes de la g\'eom\'etrie}, 
 \newblock {\em Encyclop\'edie des Sciences Math\'ematiques pures at appliqu\'ees}, Paris, 1911, Vol.\,\textbf{III.1}, pp.\,1-147.


\bibitem{gergonne} J.~Gergonne: \newblock Essai sur la th\'eorie des d\'efinitions,
\newblock {\em Annales de Math.}, \textbf{9} (1918-19), pp.\,1-36.

\newblock\url{http://archive.numdam.org/article/AMPA_1818-1819__9__1_0.pdf}

\bibitem{greco-mazzucchi} G. H. Greco, S. Mazzucchi:  
\newblock Peano on Definition of Surface Area.
\newblock ({\em forthcoming})



\bibitem{grass} H.~G. Grassmann: \newblock {\em Gesammelte Werke},  %
\newblock Teubner, Leipzig, 1894-1911.

\newblock \url{http://quod.lib.umich.edu/cgi/t/text/text-idx?c=umhistmath&idno=ABW0785}

\bibitem{grass_branch} H.~G. Grassmann: \newblock {\em A New Branch
of Mathematics, The Ausdehnungslehre of 1844 and Other Works},
\newblock Open Court, Chicago, 1995.


\bibitem{grass_E} H.~G. Grassmann: \newblock {\em Extension Theory} 
(Translation of \emph{Ausdehnungslehre}, 2nd edition, 1862), %
\newblock American Mathematical Society, 2000.

 

\bibitem{peano_vita} H.C.~Kennedy: 
\newblock {\em Life and Works of Giuseppe Peano},
\newblock D.Reidel Publ. Co, Dordrecht 1980

\newblock\url{http://www.lulu.com/items/volume_28/413000/413765/3/print/Peano_text_2.pdf}


\bibitem{lewis} A.~C.~Lewis: \newblock H.~Grassmann's 1844 Ausdehnungslehre and Schleiermacher's Dialektik, 
\newblock {\em Annals od Science}, \textbf{34} (1977), pp.\,103--162.


\bibitem{lotze_1923} A.~Lotze:
\newblock Die Grassmannsche Ausdehnungslehre, 
\newblock in  \emph{Encyklop\"adie der mathematischen Wissenschaften mit Einschluss ihrer Anwendungen}, Leipzig, 1914-1931,
Band III, T.\,1, H.\,2, pp.\,1425-1550.

\newblock\url{http://gdz.sub.uni-goettingen.de/dms/load/img/?PPN=PPN360609767&DMDID=dmdlog169}

\bibitem{maccaferri} E.~Maccaferri:  \newblock Le definizioni per astrazione e la classe di Russell,
\newblock {\em Rend. Circolo Palermo}, \textbf{35} (1913), pp.\,165--171.

\bibitem{marcolongo_1930} R.~Marcolongo: 
\newblock \emph{Analisi vettoriale generale e applicazioni}, Vol. V {\em Elettricit\`a e magnetismo},
\newblock Zanichelli, Bologna, 1930.

\bibitem{may} K.O.~May: \newblock {\em Bibliography Research Manual of the History of Mathematics},
\newblock University Toronto Press, 1973.

\bibitem{mobius} A.F.~M\"obius: \newblock {\em Der Barycentrische Calcul},
\newblock Leipzig, 1827.

\newblock\url{
http://quod.lib.umich.edu/cgi/t/text/text-idx?
c=umhistmath;idno=AAX2934.0001.001;cc=umhistmath}

\bibitem{padoa} A.~Padoa:
\newblock  Dell'astrazione matematica,
\newblock in \emph{Questioni filosofiche} (Soc. Fil. Italiana, ed.), A. F. Formiggini - Editore, Bologna, 1908, pp.\,91-104.


\bibitem{peano87} G.~Peano: \newblock {\em Applicazioni geometriche del
calcolo infinitesimale}, \newblock Fratelli Bocca, Torino, 1887.

\newblock \url{http://historical.library.cornell.edu/cgi-bin/cul.math/docviewer?did=00610002&seq=1}

\bibitem{calcolo88} G.~Peano: 
\newblock {\em Calcolo geometrico secondo Ausdehnungslehre di {H}.
  {G}rassmann}, \newblock Fratelli Bocca, Torino, 1888.


\bibitem{peano_massimi} G.~Peano: \newblock Teoremi sui massimi e minimi
geometrici, e su normali a curve e superfici, \newblock {\em Rend. Circ. Mat. Palermo}, \textbf{2} (1888), pp.\,189--192.


\newblock \url{http://www.archive.org/details/rendicontidelci15palegoog}

\bibitem{peano_area} G.~Peano: \newblock Sulla definizione dell'area di una
superficie, \newblock {\em Rend. Acc. Lincei}, \textbf{6} (1890), pp.\,54--57.


\newblock \url{http://www.archive.org/details/rendiconti02lincgoog}

\bibitem{peano1891_elementi} G.~Peano: \newblock {\em Elementi di calcolo
geometrico}, \newblock Candeletti, Torino, 1891.


\bibitem{peano1893} G.~Peano: \newblock {\em Lezioni di analisi
infintesimale}, \newblock 2 vol., Candeletti, Torino, 1893.

\bibitem{peano_quotient} G.~Peano:
\newblock  {\em Notations de logique math\'ematique (Introduction au Formulaire de 
MathŽmatiques)},
\newblock Turin, 1894.

\newblock\url{http://gdz.sub.uni-goettingen.de/dms/load/img/?IDDOC=298275}

\bibitem{peano_recensione_Grassmann} G.~Peano: 
\newblock Recensione: Hermann Grassmann's Gesammelte math. und phys. Werke,
\newblock {\em Rivista di Matematica}, \textbf{4} (1894), pp.\,167-169.




\bibitem{peano_recensione_Castellano} G.~Peano: 
\newblock Recensione: F. Castellano, Lezioni di Meccanica Razionale,
\newblock {\em Rivista di Matematica}, \textbf{5} (1895), pp.\,11-18.

\bibitem{saggio} G.~Peano: \newblock Saggio di calcolo geometrico, \newblock
{\em Atti R. Accad. Scienze Torino}, \textbf{31} (1895-96), pp.\,952--975.

\newblock\url{http://www.archive.org/details/attidellarealeac31real}


\bibitem{peano_vettori} G.~Peano: \newblock Analisi della teoria dei
vettori, \newblock {\em Atti R. Acc. Scienze Torino}, \textbf{33} (1897-98), pp.\,513--534.

\newblock\url{http://www.archive.org/details/attidellaraccad04unkngoog}


\bibitem{peano_def1900} G.~Peano: \newblock Les d\'efinitions math\'ematiques, \newblock {\em Congr\`es intern. de philosophie, Paris}, \textbf{3} (1900), pp.\,279--268.

\bibitem{peano_F2} G.~Peano:  \newblock {\em Formulaire de Math\'ematiques (tome
II)}, \newblock Bocca Fr\`eres, Paris, 1899. 

\newblock
\url{http://www.archive.org/details/formulairedemat02peangoog}

\bibitem{peano_F3} G.~Peano:  \newblock {\em Formulaire de Math\'ematiques (tome
III)}, \newblock G. Carr\`e \& C.Naud, Paris, 1901. 
\newblock
\url{http://www.archive.org/details/formulairedemat00peangoog}

\bibitem{peano_F4} G.~Peano:  \newblock {\em Formulaire Math\'ematique (tome
IV)}, \newblock Fratelli Bocca, Torino, 1903.

\newblock\url{http://www.archive.org/details/formulairedemat01peangoog}


\bibitem{peano_1908} G.~Peano: \newblock {\em Formulario Mathematico (Editio
V)}, \newblock Fratelli Bocca, Torino, 1908.

\newblock\url{http://www.archive.org/details/formulairedemat04peangoog}

\bibitem{peano_review} G.~Peano: 
\newblock  Le definizioni per astrazione,
\newblock {\em Bollettino della Mathesis}, \textbf{7} (1915), pp.\,106--120. 

\bibitem{peano_simboli} G.~Peano:
\newblock  Importanza dei simboli in matematica,
\newblock {\em Scientia, Rivista di Scienza}, \textbf{18} (1915), pp.\,165--173.

\bibitem{peano_eguale} G.~Peano: \newblock Eguale. \newblock {\em Bollettino di matematica}, \textbf{15} (1917-18), pp.\,195--198.

\bibitem{peano_def1921} G.~Peano: 
\newblock  Le definizioni in matematica, 
\newblock {\em Periodico di matematiche}, \textbf{1} (1921), pp.\,175-189.



\bibitem{peano_opere} G.~Peano: \newblock {\em Opere scelte}, 3 vol., \newblock Edizioni Cremonese, Roma 1957-9.

\bibitem{peano_1908_reprint} G.~Peano: \newblock {\em Formulario Mathematico}, Edizioni Cremonese, Roma, 1960 (reprint of \cite{peano_1908})

\bibitem{peano_english} G.~Peano:
\newblock  {\em Selected works} (H.C. Kennedy),
\newblock University Toronto Press, 1973


\bibitem{peano_omnia} G.~Peano: \newblock {\em Opera Omnia}, 1 CD-ROM (S. Boero, ed.), 
\newblock Dipartimento di Matematica, Universit\`a, Torino, 2002.

\bibitem{pensaMI} A.~Pensa:  \newblock Geometria assoluta dei vettori e delle omografie vettoriali in un $S_n$ euclideo,
\newblock {\em Rend. R. Istituto Lombardo Sc., Lett.}, \textbf{52} (1919), pp.\,339--453.

\bibitem{pensa} A.~Pensa: \newblock Geometria assoluta delle formazioni
geometriche in un $S_n$ euclideo, (I and II),
\newblock {\em Atti R. Istituto Veneto Sc., Lett. e Arti}, \textbf{79} 1919--20), pp.\,275--292, 737--761.

\bibitem{russell_principles} B.~Russell: 
\newblock  \emph{The Principles of Mathematics},
\newblock  1903.

\newblock\url{http://www.archive.org/details/principlesofmath010689mbp}

\bibitem{russell_denoting} B.~Russell: 
\newblock  On Denoting,
\newblock {\em Mind}, \textbf{14} (1905), pp.\,479-493.

\newblock\url{http://www.jstor.org/pss/2248381}

\bibitem{russell3} B.~Russell: 
\newblock  The Collected Papers, Vol. 3 (Toward the ``Principles of Mathematics'' 1900-02), 
\newblock  Routledge, London, 1993.


\bibitem{saccheri_L} G.~Saccheri:
\newblock  {\em Logica demonstrativa},
\newblock Torino, 1697 (reprinted by Olms, Hildesheim, 1980).

\bibitem{sallent} E.~Sallent Del Colombo:
\newblock  {\em CesareBurali-Forti. Contributi alla Fisica-matematica 
del primo quarto del XX secolo}, Departament de F\'\i sica,  Universitat de Barcelona,  2007. 






\bibitem{vacca} G.~Vacca: 
\newblock  Sui precursori della logica matematica I, II,
\newblock {\em Revue de Math\`ematiques}, \textbf{6} (1899), pp.\,121-125, 183-186.

\bibitem{aristotele} G~Vailati: 
\newblock  La Teoria Aristotelica della Definizione (1903), in
\newblock {\em Scritti}, Vol.\,1, pp.\,317-328.

\bibitem{saccheri} G.~Vailati: 
\newblock  Di un'opera dimenticata di P. Gerolamo Saccheri, Logica Demonstrativa (1903) in
\newblock {\em Scritti}, Vol.\,2, pp.\,212-219.

\bibitem{vailati_G} G~Vailati: 
\newblock  La grammatica dell'algebra (1908); in
\newblock {\em Scritti}, Vol.\,1, pp.\,92-110.


\bibitem{vailati} G.~Vailati:
\newblock  {\em Scritti} (M. Quaranta ed.),
\newblock Vol. 1 (Scritti di filosofia), Vol. 2 (Scritti di scienza), Vol. 3 (Scritti di scienze umane), Arnaldo Forni Editore, 1987.

\bibitem{weyl} H.~Weyl: \newblock {\em Raum-Zeit-Materie}, \newblock Springer, Berlin, 1918.
(Quoted text is extracted from the English translation: \emph{Space, Time, Matter}, Dover Publ. 1952). 



\end{thebibliography}
\end{document}